\documentclass[11pt,reqno]{amsart}

\usepackage[T1]{fontenc}
\usepackage{amsmath,amssymb,mathtools}
\usepackage{mathrsfs}
\usepackage{enumitem}
\usepackage{aliascnt}
\usepackage[colorlinks=true,linkcolor=blue,citecolor=blue,urlcolor=blue]{hyperref}
\usepackage[nameinlink,capitalise]{cleveref}

\newtheorem{theorem}{Theorem}[section]
\crefname{theorem}{Theorem}{Theorems}
\Crefname{theorem}{Theorem}{Theorems}
\newaliascnt{proposition}{theorem}
\newtheorem{proposition}[proposition]{Proposition}
\aliascntresetthe{proposition}
\crefname{proposition}{Proposition}{Propositions}
\Crefname{proposition}{Proposition}{Propositions}
\newaliascnt{lemma}{theorem}
\newtheorem{lemma}[lemma]{Lemma}
\aliascntresetthe{lemma}
\crefname{lemma}{Lemma}{Lemmas}
\Crefname{lemma}{Lemma}{Lemmas}
\newaliascnt{corollary}{theorem}
\newtheorem{corollary}[corollary]{Corollary}
\aliascntresetthe{corollary}
\crefname{corollary}{Corollary}{Corollaries}
\Crefname{corollary}{Corollary}{Corollaries}

\theoremstyle{definition}
\newaliascnt{remark}{theorem}
\newtheorem{remark}[remark]{Remark}
\aliascntresetthe{remark}
\crefname{remark}{Remark}{Remarks}
\Crefname{remark}{Remark}{Remarks}
\theoremstyle{plain}

\numberwithin{equation}{section}

\newcommand{\R}{\mathbb R}

\newcommand{\cH}{\mathcal H}
\newcommand{\cS}{\mathcal S}
\newcommand{\eps}{\varepsilon}
\newcommand{\dd}{\,\mathrm d}

\title[Normalized Mountain Pass Solution]
{Existence of a Normalized Mountain Pass Solution for a Coupled Cubic
Schr\"odinger System with Partial Confinement}

	\author[Tang, Sun]{Wan-Yue Tang$^{1}$, Hong-Rui Sun$^{1,*}$, Jianlun Liu$^{1}$ }
\thanks{\hspace{-.6cm}
	$^1$ School of Mathematics and Statistics, Lanzhou University, Lanzhou, Gansu 730000, P. R. China.
	\\
	$^{*}${\sf Corresponding author} (hrsun@lzu.edu.cn)}

\date{\today}

\begin{document}

\begin{abstract}
This paper studies a mass supercritical Schr\"odinger system with partial confinement:
\[
\begin{cases}
-\Delta u_1+|y|^2u_1
=\lambda_1u_1+\mu_1u_1^3+2\beta u_1u_2^2,\qquad x\in\R^3,\\
-\Delta u_2+|y|^2u_2
=\lambda_2u_2+\mu_2u_2^3+2\beta u_1^2u_2,\qquad x\in\R^3,
\end{cases}
\]
subject to the mass constraints \(\int_{\R^3}u_1^2\dd x=a_1\), \(\int_{\R^3}u_2^2\dd x=a_2\), where $x=(y,z)\in\R^2\times\R$, the parameters
$\mu_1,\mu_2,\beta$ are positive, and
$\lambda_1,\lambda_2\in\R$ are unknown Lagrange multipliers.
The masses satisfy $(a_1,a_2)=\eps(\rho_1,\rho_2)$ for a fixed vector $\rho=(\rho_1,\rho_2)\in (0,\infty)^2$. 
For all sufficiently small $\eps>0$, we establish the existence of a positive normalized solution $u_\eps^-$ at a mountain pass level on the full product mass sphere. This level coincides with the infimum of the energy over the negative part of the Pohozaev set. The main difficulty arises from the lack of compactness due to partial confinement. To overcome this difficulty, we develop a compactness argument that preserves the two independent mass constraints and show that every minimizing sequence for this infimum admits a strongly convergent subsequence after suitable common translations along the unconfined axial direction. We further establish the asymptotic relation $\eps E(u_\eps^-)\to m_0(\rho)$ as $\varepsilon\to0$, where
$m_0(\rho)>0$ denotes the infimum of the free energy under the corresponding Pohozaev constraint. In addition, we derive the asymptotic behavior of the kinetic and potential energies. The solution obtained here is distinct from the normalized ground states constructed by Jia et al. [Discrete Contin. Dyn. Syst., 40(2020), 2739-2766] through local minimization for sufficiently small masses.

\end{abstract}

\maketitle
\enlargethispage{2pt}

\section{Introduction}\label{sec:introduction}
Let $x=(y,z)\in\R^2\times\R$. The present paper concerns the following cooperative cubic Schr\"odinger system
\begin{equation}\label{eq:system}
\begin{cases}
-\Delta u_1+|y|^2u_1
=\lambda_1u_1+\mu_1u_1^3+2\beta u_1u_2^2,\qquad x\in\R^3,\\
-\Delta u_2+|y|^2u_2
=\lambda_2u_2+\mu_2u_2^3+2\beta u_1^2u_2,\qquad x\in\R^3,
\end{cases}
\end{equation}
subject to the two independent mass constraints
\begin{equation}\label{eq:masses}
\int_{\R^3}u_1^2\dd x=a_1,\qquad
\int_{\R^3}u_2^2\dd x=a_2.
\end{equation}
Here $\mu_1,\mu_2,\beta>0$, and $\lambda_1,\lambda_2\in\R$ are the Lagrange multipliers associated with the mass constraints. Throughout the paper, the prescribed masses satisfy
\begin{equation}\label{eq:mass-ray-intro}
(a_1,a_2)=\eps(\rho_1,\rho_2),\qquad \rho_1,\rho_2>0,
\end{equation}
where $\rho=(\rho_1,\rho_2)$ is fixed and $\eps>0$ is sufficiently small. Our goal is to establish the existence of a positive normalized solution of mountain pass type and determine its energy behavior as the masses tend to zero.

The system arises naturally from the time-dependent Schr\"odinger equations
\begin{equation}\label{eqphys}
	\begin{cases}
		-i\partial_t\Phi_1-\Delta\Phi_1+|y|^2\Phi_1
		=\mu_1|\Phi_1|^2\Phi_1+2\beta|\Phi_2|^2\Phi_1,\quad (t,x)\in\R^+\times\R^3,\\
		-i\partial_t\Phi_2-\Delta\Phi_2+|y|^2\Phi_2
		=\mu_2|\Phi_2|^2\Phi_2+2\beta|\Phi_1|^2\Phi_2,\quad (t,x)\in\R^+\times\R^3,
	\end{cases}
\end{equation}
where $t$ denotes time and $\Phi_i:\R^+\times\R^3\to\mathbb C$ are the component wave functions. Such systems describe multicomponent Bose-Einstein condensates and the propagation of mutually incoherent wave packets in nonlinear optics \cite{Frantzeskakis2010,Malomed2008}. At sufficiently low temperatures, a macroscopic fraction of bosons occupies the same quantum state. The resulting condensate dynamics can be described by coupled Gross-Pitaevskii equations under the mean-field approximation. In this setting, the coefficients $\mu_1,\mu_2$ characterize interactions within each component, whereas $\beta$ describes the interaction between the components. Under the present sign convention, positive coefficients correspond to attraction and negative coefficients to repulsion. Seeking standing waves of the form $\Phi_i(t,x)=e^{-i\lambda_i t}u_i(x)$ leads to the stationary system \eqref{eq:system}. The separate conservation of the two component masses provides the motivation for imposing the independent mass constraints \eqref{eq:masses}. The potential $|y|^2$ reflects the geometry of an elongated harmonic trap. In such a trap, confinement is much stronger in the transverse directions than along the axial direction. Neglecting the axial confinement therefore yields a potential that confines $y$ while leaving $z$ unconfined \cite{BellazziniEtAl2017}.

To formulate the stationary problem variationally, all spaces and functions are taken to be real unless otherwise specified. Unless otherwise specified, $L^p$ and $L^p_{\mathrm{loc}}$ refer to functions on $\R^3$. For two-component functions, membership and convergence in these spaces are understood componentwise. The notation $\|\cdot\|_p$ denotes the $L^p$ norm, $1\leq p\leq\infty$, using the pointwise Euclidean norm for two-component functions. The natural energy space is
$$
\cH=\Big\{u\in H^1(\R^3)\Bigm|\int_{\R^3}|y|^2|u|^2\dd x<\infty\Big\},
$$
equipped with the norm
$$
\|u\|_{\cH}^2=\int_{\R^3}\left(|\nabla u|^2+(1+|y|^2)u^2\right)\dd x.
$$
Its product space $\cH^2=\cH\times\cH$ has squared norm
$$
\|u\|_{\cH^2}^2=\|u_1\|_{\cH}^2+\|u_2\|_{\cH}^2.
$$
Given $a=(a_1,a_2)\in(0,\infty)^2$, the corresponding mass constraint set is defined by
$$
\cS_a=\left\{u=(u_1,u_2)\in\cH^2\;\middle|\;\|u_1\|_2^2=a_1,\ \|u_2\|_2^2=a_2\right\}.
$$
The energy functional involves the quantities
\begin{align*}
	K(u)&=\sum_{i=1}^2\int_{\R^3}|\nabla u_i|^2\dd x,
	&N_i(u)&=\int_{\R^3}u_i^4\dd x,\\
	C(u)&=\int_{\R^3}u_1^2u_2^2\dd x,
	& W(u)&=\sum_{i=1}^2\int_{\R^3}|y|^2u_i^2\dd x.
\end{align*}
Writing
$$
B(u)=\frac{\mu_1}{4}N_1(u)+\frac{\mu_2}{4}N_2(u)+\beta C(u),
$$
the associated energy takes the form
\begin{equation}\label{eq:energy}
E(u)=\frac12\bigl(K(u)+W(u)\bigr)-B(u).
\end{equation}
Normalized solutions of \eqref{eq:system} under \eqref{eq:masses} are precisely the critical points of $E|_{\cS_a}$. A normalized ground state is a constrained critical point whose energy is minimal among all such critical points with the same prescribed masses.

The variational structure presents two principal difficulties. The first comes from the $L^2$-supercritical nature of the cubic nonlinearities, even though they are Sobolev subcritical in dimension three. Indeed, the mass preserving scaling
\begin{equation}\label{eq:dilation-intro}
(t\star u)_i(x)=t^{\frac32}u_i(tx),\qquad i=1,2,
\end{equation}
maps $\cS_a$ into itself for every $t>0$. Along this scaling, the energy satisfies
\begin{equation}\label{eq:energy-fiber-intro}
E(t\star u)=\frac{t^2}{2}K(u)+\frac{t^{-2}}{2}W(u)-t^3B(u)\rightarrow-\infty\quad\text{as }t\to\infty.
\end{equation}
Thus, $E$ is unbounded from below on $\cS_a$, and a global minimizer cannot exist. The second difficulty concerns compactness. Under full harmonic confinement, corresponding to the potential $|y|^2+z^2$, the energy space embeds compactly into $L^p(\R^3)$ whenever $2\leq p<6$. Partial confinement changes this property because translations along the axial direction preserve the $\cH$ norm. Consequently, the embedding $\cH\hookrightarrow L^p(\R^3)$ is not compact \cite{BellazziniEtAl2017}. The construction of normalized mountain pass solutions must therefore address both the unboundedness of the constrained energy and the loss of compactness caused by axial translations. These difficulties have motivated several studies of Schr\"odinger equations with partial confinement. In the scalar setting, Bellazzini et al. \cite{BellazziniEtAl2017} established the existence of orbitally stable ground states and proposed the existence of a second critical point at a mountain pass level in \cite[Remark~1.10]{BellazziniEtAl2017}. Wei and Wu \cite{WeiWu2023} subsequently obtained a second positive solution for sufficiently small masses through a perturbation argument. Further developments concern uniqueness under strong confinement and quantitative estimates in the small-mass regime \cite{HongJin2024,SongYan2026}. A more recent contribution is due to Jeanjean and Song \cite{JeanjeanSong2025}, who established a mountain pass multiplicity result for the scalar functional
$$
E_p(u)
=
\frac12\int_{\R^3}
\left(
|\nabla u|^2+|y|^2|u|^2
\right)\dd x
-\frac1{p+1}\int_{\R^3}|u|^{p+1}\dd x,
\qquad \frac73<p<5,
$$
under a prescribed mass constraint. Their nonperturbative construction begins with solutions on expanding bounded domains. The symmetry of these solutions then provides the compactness needed to pass to the limit. This approach offers a different way to recover compactness in partially confined scalar problems. The literature on coupled systems presents additional features arising from the interaction between components. Gou \cite{Gou2018} considered the $L^2$-subcritical system
\begin{equation}\label{eq:gou-intro}
	\begin{cases}
		-\Delta u_1+|y|^2u_1
		=\lambda_1u_1+\mu_1|u_1|^{p_1-2}u_1
		+\beta r_1|u_1|^{r_1-2}u_1|u_2|^{r_2},
		\quad x\in\R^3,\\
		-\Delta u_2+|y|^2u_2
		=\lambda_2u_2+\mu_2|u_2|^{p_2-2}u_2
		+\beta r_2|u_1|^{r_1}|u_2|^{r_2-2}u_2,
		\quad x\in\R^3,
	\end{cases}
\end{equation}
where $\mu_1,\mu_2,\beta>0$, $2<p_1,p_2<\frac{10}{3}$, $r_1,r_2>1$, and $r_1+r_2<\frac{10}{3}$. In this regime, the constrained energy remains bounded from below. By exploiting this property, the author established the existence and orbital stability of normalized ground states.  Systems with quadratic interactions lead to a different setting. Alvarez and Esfahani \cite{AlvarezEsfahani2025} obtained a second normalized mountain pass solution for partially confined systems with quadratic interactions in dimensions four and five. Their construction involves a single weighted constraint on the total mass and requires a prescribed relation between the two frequencies. The quadratic interaction also has a different scaling behavior from the cubic terms within each component and the coupling terms considered here. Related existence results for partially confined Schr\"odinger-Poisson equations can be found in \cite{ShanShuaiYe2026}. Other aspects of coupled cubic systems, including invariant sets, blow-up criteria, and standing-wave instability, have been investigated in \cite{WangFeng2025,ZhouZhang2025,ZhouZhang2026}. The theory of normalized solutions for coupled systems without partial confinement has also been developed in \cite{BartschJeanjeanSoave2016,BartschZhongZou2021,ChenZou2021,GouJeanjean2018,Ikoma2014,IkomaMiyamoto2023}.

Of particular relevance to the present problem is the work of Jia et al.\ \cite{JiaLiLuo2020}. They studied the cooperative cubic system under two prescribed mass constraints and constructed constrained local minimizers when the masses are sufficiently small. These minimizers are normalized ground states, as recalled in \cref{rem:known-ground-state} below. Their result provides the starting point for the present study. While local minimization yields a positive ground state, the mountain pass structure of the full product mass sphere leads to the question of whether a second positive normalized solution exists. The present paper addresses this question by constructing a positive normalized solution at a mountain pass level associated with fixed endpoints on the full real product mass sphere. The construction also identifies this level with the infimum of the energy over the negative part of the Pohozaev set and establishes its attainment. Beyond existence, the small-mass behavior of the mountain pass energy is determined through comparison with the corresponding free system. The interaction coefficients and the vector $\rho\in(0,\infty)^2$ are arbitrary but fixed throughout the analysis; the smallness threshold may depend on these parameters.

To state our main result precisely, introduce the Pohozaev
functional $P$ and its scaling derivative $\mathcal Q$ by
\begin{align}
 P(u)&=K(u)-W(u)-3B(u),\label{eq:Pohozaev-functional}\\
 \mathcal Q(u)&=2K(u)+2W(u)-9B(u).
 \label{eq:fiber-curvature-main}
\end{align}
If $P(u)=0$, the function $t\mapsto E(t\star u)$ has a critical
point at $t=1$, where its second derivative equals $\mathcal Q(u)$.
Accordingly, we call $\{u\in\cS_a\mid P(u)=0\}$ the Pohozaev set
and define its negative part and the associated energy level by
\begin{equation}\label{eq:negative-sheet-main}
	\begin{aligned}
		&\mathcal P_a^-:=\{u\in\cS_a\mid P(u)=0,\ \mathcal Q(u)<0\},\\
		&d_\eps:=\inf_{u\in\mathcal P_{\eps\rho}^-}E(u).
	\end{aligned}
\end{equation}
The mountain pass construction uses two fixed endpoints.  Let
\begin{equation}\label{eq:gaussian-intro}
	g(x)=\pi^{-\frac34}e^{-\frac{|x|^2}{2}},
	\qquad x\in\R^3,
	\qquad \|g\|_2=1,
\end{equation}
be the $L^2$-normalized Gaussian function. Associated with each $\eps>0$ are the points
\begin{equation}\label{eq:mp-fixed-endpoints}
 \begin{split}
 p_\eps&=(\sqrt{\eps\rho_1}\,g,\sqrt{\eps\rho_2}\,g),\\
 e_\eps&=T_\eps\star p_\eps,
 \end{split}
\end{equation}
where $T_\eps=2+\frac{K(p_\eps)+W(p_\eps)}{B(p_\eps)}$. These endpoints determine the path class
\begin{equation}\label{eq:mp-level-main}
	\Gamma_\eps=\{\gamma\in C([0,1],\cS_{\eps\rho})\mid
	\gamma(0)=p_\eps,\ \gamma(1)=e_\eps\},
\end{equation}
and the associated mountain pass level is
\[
m_\eps=\inf_{\gamma\in\Gamma_\eps}
\max_{s\in[0,1]}E(\gamma(s)).
\]
For all sufficiently small $\eps>0$, the mountain pass level $m_\eps$ satisfies
$$
m_\eps>\max\{E(p_\eps),E(e_\eps)\}.
$$
The asymptotic behavior of this level is governed by a comparison problem without the confining potential. To describe it, set
$$
H=H^1(\R^3)\times H^1(\R^3)
$$
and consider the free mass constraint set
$$
\cS_\rho^0
=\left\{v=(v_1,v_2)\in H\;\middle|\;\|v_1\|_2^2=\rho_1,\ \|v_2\|_2^2=\rho_2\right\}.
$$
Using the same definitions of $K$ and $B$, the free energy and its Pohozaev functional are
\begin{equation}\label{eq:free-functionals}
	E_0(v)=\frac12K(v)-B(v),
	\qquad
	P_0(v)=K(v)-3B(v).
\end{equation}
The corresponding comparison level is defined by
\begin{equation}\label{eq:m3d}
	m_0(\rho)
	:=\inf\{E_0(v)\mid v\in\cS_\rho^0,\ P_0(v)=0\}.
\end{equation}
Moreover, it follows from the characterization \eqref{eq:free-Rayleigh} that $0<m_0(\rho)<\infty$. The comparison argument does not require the existence of a minimizer for this free problem. With these preparations, the main result can be stated as follows.

\begin{theorem}\label{thm:main}
Let $\mu_1,\mu_2,\beta>0$ and
$\rho=(\rho_1,\rho_2)\in(0,\infty)^2$ be fixed.
There exists
$\eps_0=\eps_0(\mu_1,\mu_2,\beta,\rho)>0$
such that for every $0<\eps<\eps_0$, system
\eqref{eq:system}--\eqref{eq:masses} with $a_i=\eps\rho_i$
admits a positive normalized solution
$u_\eps^-\in\mathcal P_{\eps\rho}^-$ satisfying
\begin{equation}\label{eq:main-mountain-pass-identification}
 E(u_\eps^-)=m_\eps=d_\eps.
\end{equation}
Moreover, every family of positive minimizers $u_\eps^-$ of
$d_\eps$ satisfies
\begin{equation}\label{eq:main-high-limits}
 \eps E(u_\eps^-)\to m_0(\rho),\qquad
 \eps K(u_\eps^-)\to6m_0(\rho),\qquad
 \eps W(u_\eps^-)\to0
 \quad\text{as }\eps\to0.
\end{equation}
\end{theorem}

The variational construction also leads to a compactness result for minimizing sequences on the negative part of the Pohozaev set. The lack of compactness caused by axial translations makes it natural to formulate this result modulo common translations of both components. More precisely, let
$$
\tau_s v(y,z)=v(y,z-s),\quad s\in\R,
$$
denote the common axial translation.

\begin{corollary}\label{cor:negative-sheet-compactness}
	Under the assumptions of \cref{thm:main}, fix $\eps\in(0,\eps_0)$.
	Suppose that $\{v_n\}\subset\mathcal P_{\eps\rho}^-$ is a minimizing sequence for $d_\eps$.
	Then, up to a subsequence, there exist a sequence
	$\{a_n\}\subset\R$ and a minimizer
	$V\in\mathcal P_{\eps\rho}^-$ such that
	\[
	\tau_{-a_n}v_n\to V\quad\text{in }\cH^2.
	\]
	Each component of $V$ has a constant sign.
	If both components of $v_n$ are nonnegative for every $n$,
	then both components of $V$ are strictly positive in $\R^3$.
\end{corollary}

\begin{remark}\label{rem:known-ground-state}
For comparison, recall \cite[Theorem~1.1]{JiaLiLuo2020},
whose coupling coefficient is $2\beta$ in our notation.
For every fixed $\chi>0$, the following local level has a minimizer when
$a_1+a_2>0$ is sufficiently small:
\[
 m^\chi(a)=\inf\{E(u)\mid u\in\cS_a,\ K(u)+W(u)\leq\chi\}.
\]
Every minimizer satisfies $K(u)+W(u)<\chi$ and is a
constrained critical point of $E$ on $\cS_a$.
A positive minimizer can be chosen. Moreover, for sufficiently
small masses, these minimizers are normalized ground states,
with the energy characterization
\[
 m^\chi(a)=\inf\{E(v)\mid v\in\cS_a,\ (E|_{\cS_a})'(v)=0\}.
\]
Moreover, the orbital stability of the set of local minimizers was
proved in \cite[Theorem~1.2]{JiaLiLuo2020}.

Along $a=\eps\rho$, let $u_\eps^{\mathrm{le}}$ denote a positive
normalized ground state of \eqref{eq:system}--\eqref{eq:masses}.
Its energy admits the following estimate for some $C>0$
independent of sufficiently small $\eps$:
\begin{equation}\label{eq:known-low-energy}
 0<E(u_\eps^{\mathrm{le}})\leq C\eps.
\end{equation}
Indeed, choose a fixed $\psi\in C_c^\infty(\R^3)$ with
$\|\psi\|_2=1$. For small $\eps$, the pair
$(\sqrt{\eps\rho_1}\psi,\sqrt{\eps\rho_2}\psi)$ belongs to the
local constraint and has energy $O(\eps)$.
The constrained Pohozaev identity gives $E=\frac{K+5W}{6}>0$,
which proves positivity of the critical energy.
On the other hand, \eqref{eq:main-high-limits} implies
$E(u_\eps^-)\to+\infty$ as $\eps\to0$.
Thus \cref{thm:main}, together with \cite{JiaLiLuo2020}, gives
two distinct positive normalized solutions for all sufficiently
small $\eps>0$.
\end{remark}

The proof of \cref{thm:main} consists of several steps. The starting point is the mountain pass characterization $m_\eps=d_\eps$, which is established before proving the attainment of $d_\eps$. Every admissible path crosses a boundary on which the energy is at least $d_\eps$, and this gives $m_\eps\geq d_\eps$. To obtain the reverse inequality, a minimizing sequence for $d_\eps$ is used to construct admissible paths through suitable mass preserving scalings and interpolation. The resulting upper bound completes the characterization of the mountain pass level. Once this characterization is available, the main obstacle is to establish compactness for a constrained Palais-Smale sequence at $d_\eps$, up to common axial translations. The difficulty is not limited to the translation invariance of the problem. A minimizing sequence may split into several profiles, and some nonzero profiles may have one vanishing component. The remainder may also retain mass and quadratic energy even when its quartic terms vanish (see also \cite{Ikoma2014}). Hence, the compactness argument must recover each prescribed mass and exclude any residual kinetic or potential energy.

The local minimization problem considered by Jia et al.\ \cite{JiaLiLuo2020} admits a strict subadditivity argument based on coupled rearrangements. In the present setting, the comparison pair must lie in $\mathcal P_{\eps\rho}^-$. Combining different profiles and restoring the missing masses does not necessarily preserve the Pohozaev identity. As a result, an additional estimate is needed to control the energy cost of returning to the negative part of the Pohozaev set. The componentwise combination of the nonnegative profiles through $(f_i^2+g_i^2)^{\frac12}$ does not increase the sum of their kinetic energies. Any mass not contained in these profiles is supplied by functions spread out in the axial direction, and their energy contribution is estimated separately. If splitting or mass loss occurs, this construction and a scaling with parameter $t$ close to one give
$$
t\star w\in\mathcal P_{\eps\rho}^-,\qquad
E(t\star w)\leq d_\eps-D+C_\eps D^2<d_\eps.
$$
Here $D>0$ can be made arbitrarily small, while $C_\eps>0$ is independent of $D$ for fixed $\eps$. Choosing $D$ sufficiently small gives the strict inequality, which contradicts the definition of $d_\eps$. This rules out splitting and mass loss. With these possibilities excluded, minimality further forces the remaining quadratic energy to vanish. Strong convergence in $\cH^2$ then follows after passing to a subsequence and applying suitable common axial translations. Neither the admissible paths nor the minimizing sequences require symmetry assumptions. The proof is completed by combining the energy level estimates with the Pohozaev identities to derive the small-mass limits.

The paper is organized as follows. \Cref{sec:setting} introduces the variational preliminaries. The energy estimates and the mountain pass characterization are established in \cref{sec:levels}. \Cref{sec:profiles} develops the bounded constrained Palais-Smale sequence and its finite axial profile decomposition. Finally, \cref{sec:attainment} proves the attainment of $d_\eps$, establishes the small-mass asymptotics, and completes the proof of \cref{thm:main}.

\section{Variational preliminaries}\label{sec:setting}

We use the spaces and functionals introduced in
\cref{sec:introduction}. The inner product on $\cH$ is induced by
its norm. On $\cH^2$, we use the sum of the
component inner products.

The Sobolev embedding into $L^4(\R^3)$ and H\"older's inequality give
$E\in C^1(\cH^2,\R)$, see \cite[Corollary~9.10]{Brezis2011}.
The mass functionals $M_i$ and the derivatives $DM_i(u)$ are given by
\begin{equation}\label{eq:mass-maps}
	M_i(u)=\frac12\|u_i\|_2^2,
	\qquad DM_i(u)[h]=\int_{\R^3}u_i h_i\dd x,
	\qquad i=1,2.
\end{equation}
Since $a_i>0$, the two mass gradients are linearly independent on $\cS_a$.
Hence every constrained critical point of $E|_{\cS_a}$ is a weak solution
of \eqref{eq:system} with two Lagrange multipliers.

The scaling in \eqref{eq:dilation-intro} preserves both masses and
satisfies
\begin{equation}\label{eq:scaling}
	K(t\star u)=t^2K(u),\quad W(t\star u)=t^{-2}W(u),\quad
	B(t\star u)=t^3B(u).
\end{equation}
The map $(t,u)\mapsto t\star u$ is jointly continuous from
$(0,\infty)\times\cH^2$ to $\cH^2$. Indeed, the scaling operators
are uniformly bounded when $t$ ranges over a compact subinterval of
$(0,\infty)$, while continuity for each fixed $u$ follows by approximation
with smooth compactly supported functions.
Differentiating the explicit formula \eqref{eq:energy-fiber-intro} gives
\[
P(t\star u)=t\frac{\dd}{\dd t}E(t\star u),\qquad
\mathcal Q(t\star u)=t\frac{\dd}{\dd t}P(t\star u).
\]
When $P(u)=0$, $\mathcal Q(u)$ is the second derivative of
$t\mapsto E(t\star u)$ at $t=1$. 

\begin{proposition}\label{prop:pohozaev}
	Every constrained critical point $u\in\cS_a$ of $E$ satisfies $P(u)=0$.
\end{proposition}
\begin{proof}
	A constrained critical point is a weak solution of \eqref{eq:system}.
	The constrained Pohozaev identity in
	\cite[Lemma~2.9]{JiaLiLuo2020}, with the coupling coefficient there
	taken as $2\beta$, gives
	\[
	K(u)-W(u)-\frac34\mu_1N_1(u)-\frac34\mu_2N_2(u)
	-3\beta C(u)=0.
	\]
	By the definition of $B$, the preceding identity can be written as $$P(u)=K(u)-W(u)-3B(u)=0.$$
\end{proof}

On the Pohozaev constraint $P(u)=0$,  the defining identities reduce to
\begin{equation}\label{eq:sheet-identities}
	\begin{aligned}
		& B(u) =\frac{K(u)-W(u)}3,\\
		& E(u) =\frac{K(u)+5W(u)}6,\quad
		\mathcal Q(u)=-K(u)+5W(u).
	\end{aligned}
\end{equation}
Consequently, 
$$
\mathcal P_{\eps\rho}^-=\{u\in\cS_{\eps\rho}\mid P(u)=0,~K(u)>5W(u)\}.
$$

For the transverse harmonic oscillator, integration by parts yields
the following identity for $v\in C_c^\infty(\R^2)$:
\begin{align*}
	& \int_{\R^2}\bigl(|\nabla_y v|^2+|y|^2|v|^2\bigr)\dd y =2\|v\|_{L^2(\R^2)}^2
	+\sum_{j=1}^2\|\partial_{y_j}v+y_jv\|_{L^2(\R^2)}^2\\
	& \geq 2\|v\|_{L^2(\R^2)}^2.
\end{align*}
By approximation with smooth compactly supported functions, the preceding identity holds for every \(v\in H^1(\mathbb R^2)\) satisfying
\(
\int_{\mathbb R^2}|y|^2|v(y)|^2\,\mathrm dy<\infty.
\)
Equality in the preceding inequality holds for the \(L^2(\mathbb R^2)\)-normalized Gaussian function
\begin{equation}\label{eq:transverse-ground-state}
	\phi_0(y)=\pi^{-\frac{1}{2}}\mathrm e^{-\frac{|y|^2}{2}},\qquad
	(-\Delta_y+|y|^2)\phi_0=2\phi_0.
\end{equation}
Hence the lowest eigenvalue of
$-\Delta_y+|y|^2$ on $L^2(\R^2)$ is $\Lambda_0=2$.
Integrating the transverse inequality in $z$ and using
$\|\partial_z u_1\|_2^2+\|\partial_z u_2\|_2^2\geq0$, we obtain
\begin{equation}\label{eq:spectral-bound}
	K(u)+W(u)\geq\Lambda_0(\|u_1\|_2^2+\|u_2\|_2^2)
	\qquad u\in\cH^2.
\end{equation}
The Gaussian function $g$ in \eqref{eq:gaussian-intro} is related
to $\phi_0$ by
\[
g(y,z)=\pi^{-\frac{1}{4}}\phi_0(y)\mathrm e^{-\frac{z^2}{2}}.
\]
A direct calculation gives
\[
\int_{\R^3}|\nabla g|^2\dd x=\frac32,\qquad
\int_{\R^3}|y|^2g^2\dd x=1,\qquad
\int_{\R^3}g^4\dd x=\frac1{(2\pi)^{\frac{3}{2}}}.
\]
The endpoints in \eqref{eq:mp-fixed-endpoints} are componentwise
positive and belong to $\cS_{\eps\rho}$, since $\|g\|_2=1$ and
scaling preserves both masses. These integrals yield
\[
\begin{aligned}
	& K(p_\eps)+W(p_\eps) =\frac52\eps(\rho_1+\rho_2),\\
	& B(p_\eps) =\frac{\eps^2}{(2\pi)^{\frac{3}{2}}}
	\left(\frac{\mu_1}{4}\rho_1^2+
	\frac{\mu_2}{4}\rho_2^2+\beta\rho_1\rho_2\right).
\end{aligned}
\]
In particular, $E(p_\eps)=O(\eps)$. By the definition of $T_\eps$ following \eqref{eq:mp-fixed-endpoints}, we obtain
\[
T_\eps B(p_\eps)=K(p_\eps)+W(p_\eps)+2B(p_\eps).
\]
Since $T_\eps>1$, \eqref{eq:scaling} gives
\[ \begin{aligned}
	E(e_\eps) =&\frac{T_\eps^2}{2}K(p_\eps) +\frac{T_\eps^{-2}}{2}W(p_\eps)-T_\eps^3B(p_\eps)\\
	\leq & T_\eps^2 \left(\frac{K(p_\eps)+W(p_\eps)}2-T_\eps B(p_\eps)\right)\\
	=&-T_\eps^2 \left(\frac{K(p_\eps)+W(p_\eps)}2+2B(p_\eps)\right)<0.
\end{aligned} \]
Since mass preserving scaling is continuous in \(\cH^2\), the path
\[
\gamma_\eps(s)=(1+s(T_\eps-1))\star p_\eps,
\qquad s\in[0,1],
\]
belongs to \(\Gamma_\eps\). Hence \(\Gamma_\eps\neq\varnothing\).

For the free mass sphere and functionals introduced in
\eqref{eq:free-functionals}--\eqref{eq:m3d}, define the optimal
Gagliardo-Nirenberg constant
\[
\mathcal C_\beta(\rho)
=\sup_{v\in\cS_\rho^0}
\frac{4B(v)}{K(v)^{\frac{3}{2}}}.
\]
The scalar Gagliardo-Nirenberg inequality and H\"older's inequality
show that $\mathcal C_\beta(\rho)<\infty$. Moreover, evaluating the quotient at any pair in
$\cS_\rho^0$ shows that the constant is positive. Hence
$0<\mathcal C_\beta(\rho)<\infty$, see also the corresponding
estimates in \cite[Section~1]{JiaLiLuo2020}.
For \(v\in\cS_\rho^0\), consider the energy functional \(E_0\) without the partial potential along the mass preserving scaling:
\[
\Psi_v(t)=E_0(t\star v)
=\frac{t^2}{2}K(v)-t^3B(v),\qquad t>0.
\]
By differentiation and the definition of \(P_0\) in \eqref{eq:free-functionals}, we obtain
\[
\Psi_v'(t)=tK(v)-3t^2B(v),\qquad
P_0(t\star v)=t^2K(v)-3t^3B(v)=t\Psi_v'(t).
\]
The equation $P_0(t\star v)=0$ has the unique positive solution
\[
t_v=\frac{K(v)}{3B(v)}.
\]
The function $\Psi_v$ increases on $(0,t_v)$ and decreases on
$(t_v,\infty)$. Hence $t_v$ is its global maximizer.
Evaluating $\Psi_v(t_v)$ gives
\begin{equation}\label{eq:free-Rayleigh}
	m_0(\rho)
	=\inf_{v\in\cS_\rho^0}\frac{K(v)^3}{54B(v)^2}
	=\frac{8}{27\mathcal C_\beta(\rho)^2}.
\end{equation}
In particular, $0<m_0(\rho)<\infty$. The proof does not require the existence of a minimizer for \(m_0(\rho)\).
For \(u\in\cS_{\eps\rho}\), set \(v=u/\sqrt{\eps}\in\cS_\rho^0\). By the definition of \(\mathcal C_\beta(\rho)\), we obtain
\begin{equation}\label{eq:B-GN-eps}
	B(u)\leq\frac14\eps^{\frac{1}{2}}\mathcal C_\beta(\rho)K(u)^{\frac{3}{2}}
	\qquad(u\in\cS_{\eps\rho}).
\end{equation}

\section{Level estimates and the mountain pass characterization}
\label{sec:levels}

In this section, we estimate $d_\eps$ and prove its mountain pass characterization.

\begin{lemma}\label{lem:fold-gap}
	Fix $0<\vartheta<\widehat\vartheta<\frac{7}{25}$. There exists
	$\eps_0>0$, depending only on $\mu_1,\mu_2,\beta,\rho,\vartheta$,
	such that for every $\eps\in(0,\eps_0)$, the set
	$\mathcal P_{\eps\rho}^-$ is nonempty and
	\begin{equation}\label{eq:d-upper-rough}
		\frac{m_0(\rho)}\eps\leq d_\eps
		\leq\frac{1+\vartheta}{\eps}m_0(\rho).
	\end{equation}
	Moreover, there exist constants \(r_*\in(0,\frac15)\) and \(\kappa>0\) with the following property. For every \(\eps\in(0,\eps_0)\) and \(u\in\mathcal P_{\eps\rho}^-\), if
	\begin{equation}\label{eq:near-d-fold}
		E(u)\leq\frac{1+\widehat\vartheta}{\eps}m_0(\rho)
	\end{equation}
	then
	\begin{equation}\label{eq:uniform-negative-curvature}
		\frac{W(u)}{K(u)}\leq r_*<\frac{1}{5},
		\qquad \mathcal Q(u)\leq-\frac\kappa\eps.
	\end{equation}
	The constant $r_*$ depends only on $\widehat\vartheta$,
	whereas $\kappa$ may also depend on the fixed interaction
	parameters and the vector $\rho$.
\end{lemma}

\begin{proof}
	First, we prove that $\mathcal P_{\eps\rho}^-$ is nonempty.
	To this end, we choose a near minimizer $w\in\cS_\rho^0$ of the quotient
	in \eqref{eq:free-Rayleigh}.
	By smooth approximation in $H^1$ and normalization of each component
	to have squared $L^2$ norm $\rho_i$, we may take
	$w\in C_c^\infty(\R^3;\R^2)\cap\cS_\rho^0$.
	Set $t_w=\frac{K(w)}{3B(w)}$ and $v=t_w\star w$.
	Then $E_0(v)=\frac{K(w)^3}{54B(w)^2}$, so a suitable choice of $w$ gives
	$v\in C_c^\infty(\R^3;\R^2)\cap\cS_\rho^0$ satisfying
	\begin{equation}\label{eq:free-near-minimizer-high}
		P_0(v)=0,\qquad
		E_0(v)<(1+\frac{\vartheta}{2})m_0(\rho).
	\end{equation}
	For $\eps>0$, define
	\(
	U_\eps(x)=\eps^{-1}v(\frac{x}{\eps}),
	\)
	where $v$ is fixed independently of $\eps$.
	By the change of variables $x=\eps\xi$, we obtain
	\begin{align*}
		& \|U_{i,\eps}\|_2^2 =\eps\rho_i, & & K(U_\eps) =\eps^{-1}K(v),\\
		& B(U_\eps) =\eps^{-1}B(v), & & W(U_\eps) =\eps^3W(v).
	\end{align*}
	Since $K(v)=3B(v)$, the equation $P(t\star U_\eps)=0$ is equivalent to
	\begin{equation}\label{eq:trial-projection}
		t^4(1-t)=\eps^4\frac{W(v)}{K(v)}.
	\end{equation}
	The function $t\mapsto t^4(1-t)$ is strictly decreasing on
	$[\frac{9}{10},1]$ and vanishes at $t=1$.
	Thus, for sufficiently small $\eps>0$, \eqref{eq:trial-projection}
	has a unique solution $t_\eps\in(\frac{9}{10},1)$, with
	$t_\eps\to1$ as $\eps\to0$.
	Then equation \eqref{eq:trial-projection} gives
	\[
	\eps\mathcal Q(t_\eps\star U_\eps)
	=t_\eps^2K(v)(4-5t_\eps)<0.
	\]
	Hence $t_\eps\star U_\eps\in\mathcal P_{\eps\rho}^-$, which proves
	that $\mathcal P_{\eps\rho}^-\neq\varnothing$.
	
	Next, we prove the upper bound for $d_\eps$.
	Since $P_0(v)=0$, the preceding scaling calculation gives
	\[
	E_0(t_\eps\star v)\leq E_0(v).
	\]
	Therefore,
	\begin{equation}\label{eq:trial-energy-bound}
		E(t_\eps\star U_\eps)
		=\eps^{-1}E_0(t_\eps\star v)
		+\frac{\eps^3}{2t_\eps^2}W(v)
		\leq\eps^{-1}E_0(v)+\frac{50}{81}\eps^3W(v).
	\end{equation}
	For sufficiently small $\eps$, we have
	$\frac{50}{81}\eps^4W(v)\leq\frac{\vartheta}{2}m_0(\rho)$.
	Combining \eqref{eq:free-near-minimizer-high} with
	\eqref{eq:trial-energy-bound} yields
	\[
	d_\eps\leq E(t_\eps\star U_\eps)
	\leq\frac{1+\vartheta}{\eps}m_0(\rho).
	\]
	This is the upper bound in \eqref{eq:d-upper-rough}.
	
	For the lower bound, let $u\in\mathcal P_{\eps\rho}^-$ and write
	$r=\frac{W(u)}{K(u)}\in[0,\frac{1}{5})$. By \eqref{eq:sheet-identities} and
	\eqref{eq:B-GN-eps},
	\[
	\frac{1-r}{3}K(u)
	\leq\frac14\eps^{\frac{1}{2}}\mathcal C_\beta(\rho)K(u)^{\frac{3}{2}},
	\qquad
	K(u)\geq\frac{16(1-r)^2}{9\eps\mathcal C_\beta(\rho)^2}.
	\]
	Substitution in the energy identity gives
	\begin{equation}\label{eq:fiber-lower-f}
		E(u)\geq\frac{m_0(\rho)}\eps f(r),
		\qquad f(r)=(1-r)^2(1+5r).
	\end{equation}
	The polynomial satisfies
	\begin{align*}
		& f(0)=1,\qquad f(\frac{1}{5})=1+\frac{7}{25},\\
		& f'(r)=3(1-r)(1-5r)>0\quad(0<r<\frac{1}{5}).
	\end{align*}
	Since $f(r)\geq1$ on $[0,\frac15]$, taking the infimum over
	$u\in\mathcal P_{\eps\rho}^-$ gives the lower bound in
	\eqref{eq:d-upper-rough}.
	
	Finally, we prove \eqref{eq:uniform-negative-curvature}.
	If $u$ satisfies \eqref{eq:near-d-fold}, then
	$f(r)\leq1+\widehat\vartheta$. Let $r_*\in(0,\frac{1}{5})$ be the
	unique solution of $f(r_*)=1+\widehat\vartheta$. It follows that
	$r\leq r_*$ and
	\[
	\mathcal Q(u)=-(1-5r)K(u)
	\leq-\frac\kappa\eps,\qquad
	\kappa=\frac{16(1-5r_*)(1-r_*)^2}{9\mathcal C_\beta(\rho)^2}>0.
	\]
	Thus both inequalities in \eqref{eq:uniform-negative-curvature} hold.
\end{proof}

For the rest of the paper, we fix
\[
\vartheta=\frac1{20},\qquad \widehat\vartheta=\frac1{10},
\]
and use the corresponding constants from \cref{lem:fold-gap}.

Next, we construct a continuous path between two nonnegative pairs
on the same product mass sphere and estimate its energy.

\begin{lemma}
	\label{lem:quadrature-path-mp}
	Let $a_1,a_2>0$ and let $f,g\in\cS_a$ be componentwise nonnegative.
	For $s\in[0,1]$, define
	\[
	(\mathcal I_s(f,g))_i
	:=\bigl((1-s)f_i^2+sg_i^2\bigr)^{\frac{1}{2}},\qquad i=1,2.
	\]
	Then $s\mapsto\mathcal I_s(f,g)$ is a continuous path from $f$ to $g$
	in $\cS_a$ with respect to the $\cH^2$ norm.
	For every $s\in[0,1]$, the following relations hold:
	\begin{align}
		& K(\mathcal I_s(f,g)) \leq(1-s)K(f)+sK(g),\label{eq:quadrature-path-kinetic}\\
		& W(\mathcal I_s(f,g)) =(1-s)W(f)+sW(g),\label{eq:quadrature-path-quadratic}\\
		& B(\mathcal I_s(f,g)) \geq(1-s)^2B(f)+s^2B(g).
		\label{eq:quadrature-path-quartic}
	\end{align}
\end{lemma}

\begin{proof}
	Write $h(s)=\mathcal I_s(f,g)$. Since
	$h_i(s)^2=(1-s)f_i^2+sg_i^2$, the component masses satisfy
	\[
	\|h_i(s)\|_2^2=(1-s)\|f_i\|_2^2+s\|g_i\|_2^2=a_i.
	\]
	Thus $h(s)\in\cS_a$, while the componentwise nonnegativity of $f$ and
	$g$ gives $h(0)=f$ and $h(1)=g$.
	
	For $0<s<1$, the chain rule for weak derivatives gives the following
	formula almost everywhere on $\{h_i(s)>0\}$:
	\[
	\nabla h_i(s)
	=\frac{(1-s)f_i\nabla f_i+sg_i\nabla g_i}
	{\bigl((1-s)f_i^2+sg_i^2\bigr)^{\frac{1}{2}}}.
	\]
	On $\{f_i=g_i=0\}$,
	the weak gradients of $f_i$, $g_i$, and $h_i(s)$ vanish almost
	everywhere. The Cauchy-Schwarz inequality gives
	\[
	|\nabla h_i(s)|^2
	\leq(1-s)|\nabla f_i|^2+s|\nabla g_i|^2.
	\]
	Integrating this inequality over $\R^3$ and summing over
	$i=1,2$ gives \eqref{eq:quadrature-path-kinetic}.
	Moreover, the identity for $h_i(s)^2$ gives
	\[
	W(h(s))=(1-s)W(f)+sW(g).
	\]
	Expanding the quartic terms in the definition of $B$, we obtain
	\[
	\begin{aligned}
		&B(h(s))-(1-s)^2B(f)-s^2B(g)\\
		=&s(1-s)\sum_{i=1}^2\frac{\mu_i}{2}
		\int_{\R^3}f_i^2g_i^2\dd x +\beta s(1-s)\int_{\R^3}
		\bigl(f_1^2g_2^2+g_1^2f_2^2\bigr)\dd x.
	\end{aligned}
	\]
	The right hand side is nonnegative since $0\leq s\leq1$ and
	$\mu_1,\mu_2,\beta>0$. This proves \eqref{eq:quadrature-path-quartic}.
	
	Finally, it remains to verify continuity in $\cH^2$.  For $s,s_0\in[0,1]$,
	\[
	|h_i(s)-h_i(s_0)|^2
	\leq |s-s_0|(f_i^2+g_i^2).
	\]
	After integration, this estimate proves continuity in both $L^2$ and the weighted space
	$L^2(\R^3,|y|^2\dd x)$. The gradient formula is pointwise continuous
	in $s\in(0,1)$ almost everywhere. Since $\nabla f_i=0$ almost
	everywhere on $\{f_i=0\}$, as $s\to0^+$ we have
	\[
	\begin{aligned}
		& \nabla h_i(s)\longrightarrow\nabla f_i
		&&\text{a.e. on }\{f_i>0\},\\
		& \nabla h_i(s)=\sqrt{s}\,\nabla g_i\longrightarrow0=\nabla f_i
		&&\text{a.e. on }\{f_i=0\}.
	\end{aligned}
	\]
	The limit as $s\to1^-$ follows similarly.
	Moreover,
	\[
	|\nabla h_i(s)|^2\leq|\nabla f_i|^2+|\nabla g_i|^2.
	\]
	Therefore, applying dominated convergence to the squared gradient
	differences gives continuity of the gradients in $L^2$.
	Together with the preceding $L^2$ and weighted $L^2$ continuity,
	this proves continuity of $s\mapsto h(s)$ in $\cH^2$.
\end{proof}

For $p,e\in\cS_{\eps\rho}$, define
\begin{align*}
 & \Gamma_\eps(p,e):=\{\gamma\in C([0,1],\cS_{\eps\rho})\mid
 \gamma(0)=p,\ \gamma(1)=e\},\\
 & m_\eps(p,e):=\inf_{\gamma\in\Gamma_\eps(p,e)}
 \max_{s\in[0,1]}E(\gamma(s)).
\end{align*}
Here path continuity is understood in the $\cH^2$ norm.

\begin{proposition}
	\label{prop:mp-level-identification}
	Let $\mu_1,\mu_2,\beta>0$ and $\rho\in(0,\infty)^2$ be fixed.
	For sufficiently small $\eps>0$, let $p,e\in\cS_{\eps\rho}$ be
	componentwise nonnegative and satisfy
	\begin{equation}\label{eq:mp-admissible-endpoints}
		B(p)<\frac4{15}K(p),\qquad
		\frac{K(p)+W(p)}2<\frac{m_0(\rho)}\eps,
		\qquad E(e)<0.
	\end{equation}
	Then $\Gamma_\eps(p,e)\neq\varnothing$ and $m_\eps(p,e)\in\R$.
	Moreover,
	\begin{equation}\label{eq:mp-equals-negative-sheet}
		m_\eps(p,e)=d_\eps>\max\{E(p),E(e)\}.
	\end{equation}
	In particular, $m_\eps=d_\eps$.
\end{proposition}

\begin{proof}
	By \cref{lem:quadrature-path-mp},
	\[
	\gamma_0(s):=\mathcal I_s(p,e),\qquad s\in[0,1],\qquad
	\gamma_0\in\Gamma_\eps(p,e).
	\]
	Hence $\Gamma_\eps(p,e)\neq\varnothing$.
	Since $E$ is continuous, its values along $\gamma_0$ are bounded.
	Together with the endpoint bound $m_\eps(p,e)\geq\max\{E(p),E(e)\}$,
	this implies $m_\eps(p,e)\in\R$.
	In this proof, write $m_0:=m_0(\rho)$ and
	$\mathcal C_\beta:=\mathcal C_\beta(\rho)$.
	Since $\vartheta=\frac1{20}$, the upper bound in \cref{lem:fold-gap}
	has coefficient $1+\vartheta=\frac{21}{20}$. Hence
	\begin{equation}\label{eq:mp-rough-levels}
		\frac{m_0}\eps\leq d_\eps
		\leq\frac{21}{20}\frac{m_0}\eps
		<\frac{32}{25}\frac{m_0}\eps.
	\end{equation}
	
	\smallskip
	\noindent\textbf{Step 1.}
	First, we prove $m_\eps(p,e)\geq d_\eps$.
	When $P(u)=0$, we have $W(u)=K(u)-3B(u)$. Therefore,
	\[
	\mathcal Q(u)=2K(u)+2W(u)-9B(u)=4K(u)-15B(u).
	\]
	When the Pohozaev constraint holds, $\mathcal Q(u)<0$ is equivalent to
	$B(u)>\frac{4K(u)}{15}$. This explains the threshold used below.
	Consider the relatively open subset of $\cS_{\eps\rho}$
	\[
	\mathcal A_\eps
	:=\left\{B<\frac4{15}K\right\}\cup\{P>0\}.
	\]
	The complement of $\mathcal A_\eps$ in $\cS_{\eps\rho}$ is
	\[
	\mathcal A_\eps^c
	=\left\{B\geq\frac4{15}K\right\}\cap\{P\leq0\}.
	\]
	Every point on the relative boundary belongs to this complement. At
	least one of the two defining inequalities must be an equality.
	Otherwise, the point would lie in the interior of $\mathcal A_\eps^c$.
	Consequently, the relative boundary satisfies
	\begin{equation}\label{eq:mp-boundary-inclusion}
		\partial_{\cS_{\eps\rho}}\mathcal A_\eps
		\subset
		\left\{P=0,\ B>\frac4{15}K\right\}
		\cup
		\left\{B=\frac4{15}K,\ P\leq0\right\}.
	\end{equation}
	For a point in the first set in \eqref{eq:mp-boundary-inclusion},
	$W=K-3B<\frac{K}{5}$. Hence $\mathcal Q=-K+5W<0$.  Such a point belongs to
	$\mathcal P_{\eps\rho}^-$ and has energy at least $d_\eps$.
	On the second set, $W\geq \frac{K}{5}$, and
	\[
	E=\frac7{30}K+\frac12W\geq\frac13K.
	\]
	By the Gagliardo-Nirenberg inequality \eqref{eq:B-GN-eps}, we obtain
	\[
	\frac4{15}K\leq\frac14\eps^{\frac{1}{2}}\mathcal C_\beta K^{\frac{3}{2}},
	\qquad K\geq\frac{256}{225\eps \mathcal C_\beta^2}.
	\]
	Using $m_0=\frac{8}{27\mathcal C_\beta^2}$, we conclude that
	\begin{equation}\label{eq:mp-artificial-boundary-gap}
		E\geq\frac{256}{675\eps \mathcal C_\beta^2}
		=\frac{32}{25}\frac{m_0}\eps>d_\eps.
	\end{equation}
	Thus every point of $\partial_{\cS_{\eps\rho}}\mathcal A_\eps$ has energy at
	least $d_\eps$. The first endpoint condition places $p$ in $\mathcal A_\eps$.
	Due to $E(e)<0$, we have
	\[
	B(e)>\frac12(K(e)+W(e))>\frac4{15}K(e),\qquad
	P(e)<-\frac12K(e)-\frac52W(e)<0.
	\]
	These strict inequalities show that $e$ is an interior point of
	$\mathcal A_\eps^c$. Since $p\in\mathcal A_\eps$, every path in
	$\Gamma_\eps(p,e)$ intersects the relative boundary.
	Consequently, the boundary energy estimate gives
	\begin{equation}\label{eq:mp-lower-comparison}
		m_\eps(p,e)\geq d_\eps.
	\end{equation}
	
	\smallskip
	\noindent\textbf{Step 2.}
	Let $\{u_n\}\subset\mathcal P_{\eps\rho}^-$ be a minimizing sequence
	for $d_\eps$. For all sufficiently large $n$, we construct a continuous
	path $\alpha_n$ from $p$ to $\sigma\star u_n$ such that
	\[
	\max_{s\in[0,1]}E(\alpha_n(s))<E(u_n),
	\]
	where $\sigma>0$ will be chosen below.
	Replacing each \(u_{i,n}\) by \(|u_{i,n}|\) does not change the prescribed masses or the values of \(K,W,B\). Consequently, the resulting sequence is still a minimizing sequence for \(d_\eps\) on \(\mathcal P_{\eps\rho}^-\), and we may assume \(u_{i,n}\geq0\).
	Since $E(u_n)\to d_\eps$ and $d_\eps\leq\frac{21}{20}\frac{m_0}{\eps}$,
	for all sufficiently large $n$ we have
	\[
	E(u_n)\leq\frac{11}{10}\frac{m_0}\eps.
	\]
	Set $r_n=\frac{W(u_n)}{K(u_n)}$.  The lower bound
	\eqref{eq:fiber-lower-f} gives
	\[
	f(r_n)\leq\frac{11}{10},\qquad
	f(r)=(1-r)^2(1+5r).
	\]
	Since $f$ is strictly increasing on $[0,\frac{1}{5}]$ and
	\[
	f(0)=1<\frac{11}{10}<f(\frac{1}{25})=\frac{3456}{3125},
	\]
	there is a unique $\bar r\in(0,\frac{1}{25})$ satisfying
	$f(\bar r)=\frac{11}{10}$. This choice is independent of $n$ and $\eps$,
	and $r_n\leq\bar r$ for all sufficiently large $n$.
	Fix $\sigma=\frac{1}{\sqrt{5}}$ and put $q_n=\sigma\star u_n$. This choice gives
	\[
	K(q_n)=\frac15K(u_n),\qquad W(q_n)=5W(u_n).
	\]
	Combining these relations with \eqref{eq:sheet-identities} gives the strict
	estimate
	\begin{equation}\label{eq:mp-low-connection-margin}
		\begin{aligned}
			& E(u_n)-\frac12\bigl(K(q_n)+W(q_n)\bigr)\\
			=&\frac{K(u_n)+5W(u_n)}6
			-\frac12\left(\frac15K(u_n)+5W(u_n)\right)\\
			=&\frac{K(u_n)-25W(u_n)}{15}>0.
		\end{aligned}
	\end{equation}
	By \cref{lem:quadrature-path-mp}, the path
	\[
	\alpha_n(s):=\mathcal I_s(p,q_n),\qquad s\in[0,1],
	\]
	is continuous in $\cH^2$ and satisfies
	\[
	\alpha_n(0)=p,\qquad \alpha_n(1)=q_n.
	\]
	Using \eqref{eq:mp-admissible-endpoints}, \eqref{eq:mp-rough-levels},
	and the definition of $d_\eps$, we obtain
	\[
	\frac12\bigl(K(p)+W(p)\bigr)
	<\frac{m_0}{\eps}\leq d_\eps\leq E(u_n),
	\]
	where the last inequality follows from
	$u_n\in\mathcal P_{\eps\rho}^-$ and
	$d_\eps=\inf_{v\in\mathcal P_{\eps\rho}^-}E(v)$.  For the other endpoint,
	\eqref{eq:mp-low-connection-margin} gives
	\[
	\frac12\bigl(K(q_n)+W(q_n)\bigr)<E(u_n).
	\]
	Now, since $B\geq0$, the definition of $E$ and the quadratic interpolation
	estimates in \cref{lem:quadrature-path-mp} yield, for every
	$s\in[0,1]$,
	\begin{align*}
		E(\alpha_n(s)) \leq & \frac12\bigl(K(\alpha_n(s))+W(\alpha_n(s))\bigr)\\
		\leq&(1-s)\frac12\bigl(K(p)+W(p)\bigr)
		+s\frac12\bigl(K(q_n)+W(q_n)\bigr)\\
		<&E(u_n).
	\end{align*}
	\smallskip
	\noindent\textbf{Step 3.}
	We show that $E(t\star u_n)\leq E(u_n)$ for every $t\geq\sigma$.
	For $u=u_n$, write $K=K(u)$ and $r=\frac{W(u)}{K(u)}$.
	Since $P(u)=0$,
	\[
	\Phi_u(t):=E(t\star u)
	=K\left(\frac12t^2+\frac r2t^{-2}-\frac{1-r}{3}t^3\right).
	\]
	Differentiating $\Phi_u(t)$, we obtain
	\begin{align*}
		& \Phi_u'(t)=Kt^{-3}h_r(t),\qquad
		h_r(t)=t^4-(1-r)t^5-r,\\
		& h_r'(t)=t^3\bigl(4-5(1-r)t\bigr).
	\end{align*}
	Since $W(u)>0$, we have $0<r<\frac1{25}$.
	Set $t_0=\frac{4}{5(1-r)}<1$. The function $h_r$ is strictly
	increasing on $(0,t_0)$ and strictly decreasing on $(t_0,\infty)$.
	Since $h_r(0)=-r<0<h_r(t_0)$ and $h_r(1)=0$, this monotonicity
	shows that $h_r(t)=0$ has exactly two positive solutions:
	$\tau_r\in(0,t_0)$ and $1$. Consequently,
	\[
	\begin{aligned}
		& \Phi_u'(t) <0, & & t\in(0,\tau_r)\cup(1,\infty),\\
		& \Phi_u'(t) >0, & & t\in(\tau_r,1).
	\end{aligned}
	\]
	Thus the maximum of $\Phi_u$ on $[\sigma,1]$ is attained
	at an endpoint. By \eqref{eq:mp-low-connection-margin},
	\[
	\Phi_u(\sigma)=E(q_n)
	\leq\frac12\bigl(K(q_n)+W(q_n)\bigr)
	<E(u_n)=\Phi_u(1).
	\]
	Since $\Phi_u$ decreases on $(1,\infty)$, we obtain
	\begin{equation}\label{eq:mp-truncated-fiber-maximum}
		\max_{t\in[\sigma,R]}E(t\star u_n)=E(u_n)
		\qquad\text{for every }R\geq1.
	\end{equation}
	Since $E(t\star u_n)\longrightarrow+\infty$ as $t\to0^+$, we take
	$t\in[\sigma,R]$ for the scaling part of the path.
	
	\smallskip
	\noindent\textbf{Step 4.}
	We connect $R_n\star u_n$ to $e$ by two paths along which the energy
	is negative. The value of $R_n\geq1$ will be chosen below.
	Set
	\begin{align*}
		& h_n(s):=\mathcal I_s(u_n,e),\qquad s\in[0,1],\\
		& L_n:=\max\bigl\{K(u_n)+W(u_n),K(e)+W(e)\bigr\},\\
		& c_n:=\frac12\min\{B(u_n),B(e)\}>0.
	\end{align*}
	By \cref{lem:quadrature-path-mp} and $(1-s)^2+s^2\geq\frac12$,
	for every $s\in[0,1]$ we have
	\[
	\begin{aligned}
		& K(h_n(s))+W(h_n(s))\leq L_n,\\
		& B(h_n(s))\geq(1-s)^2B(u_n)+s^2B(e)\\
		& \hphantom{B(h_n(s))}
		\geq\bigl((1-s)^2+s^2\bigr)\min\{B(u_n),B(e)\}\\
		& \hphantom{B(h_n(s))}\geq c_n.
	\end{aligned}
	\]
	Choose
	\[
	R_n>\max\left\{1,\frac{L_n}{2c_n}\right\}.
	\]
	For every $s\in[0,1]$, we use $R_n^{-2}\leq R_n^2$ to obtain
	\[
	\begin{aligned}
		E(R_n\star h_n(s)) &=\frac12\bigl(R_n^2K(h_n(s))+R_n^{-2}W(h_n(s))\bigr)
		-R_n^3B(h_n(s))\\
		&\leq  R_n^2\left(\frac{L_n}{2}-R_nc_n\right)<0.
	\end{aligned}
	\]
	Hence $s\mapsto R_n\star h_n(s)$ is a continuous path
	from $R_n\star u_n$ to $R_n\star e$ with $E<0$.
	For every $t\in[1,R_n]$,
	\[
	\frac{E(t\star e)}{t^3}
	=\frac{K(e)}{2t}+\frac{W(e)}{2t^5}-B(e)
	\leq\frac12\bigl(K(e)+W(e)\bigr)-B(e)=E(e)<0.
	\]
	Therefore, decreasing $t$ from $R_n$ to $1$ gives a continuous
	path from $R_n\star e$ to $e$ with $E<0$.
	
	Joining the four paths after linear reparametrization gives
	$\gamma_n\in\Gamma_\eps(p,e)$:
	\[
	p\longrightarrow\sigma\star u_n
	\longrightarrow R_n\star u_n
	\longrightarrow R_n\star e
	\longrightarrow e.
	\]
	Steps~2--4 give
	\[
	m_\eps(p,e)
	\leq\max_{s\in[0,1]}E(\gamma_n(s))
	=E(u_n).
	\]
	Since $E(u_n)\to d_\eps$, letting $n\to\infty$ yields
	\[
	m_\eps(p,e)\leq d_\eps.
	\]
	Together with \eqref{eq:mp-lower-comparison}, we obtain
	\[
	m_\eps(p,e)=d_\eps.
	\]
	Moreover, the endpoint assumptions and \eqref{eq:mp-rough-levels} give
	\[
	E(p)\leq\frac12\bigl(K(p)+W(p)\bigr)
	<\frac{m_0}{\eps}\leq d_\eps,
	\qquad E(e)<0<d_\eps.
	\]
	This proves \eqref{eq:mp-equals-negative-sheet}.
	
	For the Gaussian endpoints in \eqref{eq:mp-fixed-endpoints}, we have
	\[
	K(p_\eps)+W(p_\eps)=O(\eps),\qquad
	\frac{B(p_\eps)}{K(p_\eps)}=O(\eps).
	\]
	Together with $m_0>0$ and $E(e_\eps)<0$, these estimates verify
	\eqref{eq:mp-admissible-endpoints} for all sufficiently small $\eps>0$.
	Hence $m_\eps=m_\eps(p_\eps,e_\eps)=d_\eps$.
\end{proof}

\begin{corollary}\label{cor:mp-from-low-solution}
	For sufficiently small $\eps>0$, let
	$u_\eps^{\mathrm{le}}\in\cS_{\eps\rho}$ be a positive normalized
	ground state given by \cite[Theorem~1.1]{JiaLiLuo2020}.
	Then, for every componentwise nonnegative $e\in\cS_{\eps\rho}$
	with $E(e)<0$, we have
	\[
	m_\eps(u_\eps^{\mathrm{le}},e)=d_\eps.
	\]
\end{corollary}

\begin{proof}
	First, $P(u_\eps^{\mathrm{le}})=0$ follows from \cref{prop:pohozaev}.
	Using \eqref{eq:sheet-identities} and \eqref{eq:known-low-energy},
	we obtain
	\[
	K(u_\eps^{\mathrm{le}})+W(u_\eps^{\mathrm{le}}) \leq K(u_\eps^{\mathrm{le}})+5W(u_\eps^{\mathrm{le}})
	=6E(u_\eps^{\mathrm{le}})\leq6C\eps.
	\]
	Then the Gagliardo-Nirenberg inequality \eqref{eq:B-GN-eps} gives
	\[
	\frac{B(u_\eps^{\mathrm{le}})}{K(u_\eps^{\mathrm{le}})}
	\leq\frac14\mathcal C_\beta(\rho)
	\bigl(\eps K(u_\eps^{\mathrm{le}})\bigr)^{\frac12}
	=O(\eps).
	\]
	Since $m_0(\rho)>0$ is independent of $\eps$, these estimates imply
	that for all sufficiently small $\eps>0$,
	\[
	\begin{aligned}
		& B(u_\eps^{\mathrm{le}})<\frac4{15}K(u_\eps^{\mathrm{le}}),\\
		& \frac12\bigl(K(u_\eps^{\mathrm{le}})+W(u_\eps^{\mathrm{le}})\bigr)
		\leq3C\eps<\frac{m_0(\rho)}{\eps}.
	\end{aligned}
	\]
	Thus $u_\eps^{\mathrm{le}}$ satisfies the first two inequalities
	in \eqref{eq:mp-admissible-endpoints}.
	For every componentwise nonnegative $e\in\cS_{\eps\rho}$ with $E(e)<0$,
	\cref{prop:mp-level-identification} applies and gives the conclusion.
\end{proof}

\section{Palais-Smale sequences and axial profiles}\label{sec:profiles}

First, we construct a bounded Palais-Smale sequence for
$E|_{\cS_{\eps\rho}}$ at the level $d_\eps$.
The assumption $u\in\cH^2$ does not ensure that the formal scaling
derivative $\frac{3}{2}u+x\cdot\nabla u$ at $t=1$ belongs to $\cH^2$.
Therefore, following the scaling approach in
\cite{Jeanjean1997,BartschJeanjeanSoave2016}, we treat $t$ as an
independent variable.  On $(0,\infty)\times\cH^2$, define
\begin{align}
	& \widetilde E(t,u) :=E(t\star u)
	=\frac12\bigl(t^2K(u)+t^{-2}W(u)\bigr)-t^3B(u),
	\label{eq:augmented-energy}\\
	& \widetilde P(t,u) :=P(t\star u)=t\partial_t\widetilde E(t,u),\notag\\
	& \widetilde{\mathcal Q}(t,u) :=\mathcal Q(t\star u)=t\partial_t\widetilde P(t,u).
	\label{eq:augmented-PQ}
\end{align}
The functionals in \eqref{eq:augmented-energy}--\eqref{eq:augmented-PQ}
are $C^1$ in $(t,u)$.

Here $e_1=(1,0)$ and $e_2=(0,1)$ are the standard basis vectors of
$\R^2$. For a scalar function $f$, write $(f)_-:=\max\{-f,0\}$
for its negative part.

\begin{proposition}
	\label{prop:high-PS}
	Let $\eps\in(0,\eps_0)$. Then there exist a sequence
	$\{u_n\}\subset\mathcal P_{\eps\rho}^-$ bounded in $\cH^2$ and
	bounded sequences of
	Lagrange multipliers $\{\lambda_{i,n}\}\subset\R$, $i=1,2$, such that
	\begin{align}
		& E(u_n) \longrightarrow d_\eps\quad\text{as }n\to\infty, & & \mathcal Q(u_n) \leq-\frac\kappa{2\eps},
		\label{eq:high-PS-level}
	\end{align}
	\begin{multline}
		E'(u_n)-\lambda_{1,n}u_{1,n}e_1
		-\lambda_{2,n}u_{2,n}e_2\longrightarrow0
		\quad\text{in }(\cH^2)'
		\quad\text{as }n\to\infty.
		\label{eq:true-product-PS}
	\end{multline}
	Moreover, the negative parts satisfy
	\begin{equation}\label{eq:high-PS-asymptotic-nonnegative}
		\sum_{i=1}^2\|(u_{i,n})_-\|_2^2\longrightarrow0.
	\end{equation}
\end{proposition}

\begin{proof}
	Let $v_n\in\mathcal P_{\eps\rho}^-$ be a minimizing sequence.
	Replacing each $v_{i,n}$ by $|v_{i,n}|$ does not change the prescribed
	masses or the values of $K,W,B$. Consequently, the resulting sequence
	is still a minimizing sequence for $d_\eps$ on $\mathcal P_{\eps\rho}^-$,
	and we may assume $v_{i,n}\geq0$ before applying Ekeland's principle.
	Set $\delta_n=E(v_n)-d_\eps+\frac{1}{n}$.
	For all sufficiently large $n$, \cref{lem:fold-gap} gives
	\[
	E(v_n)\leq\frac{1+\widehat\vartheta}{\eps}m_0(\rho),
	\qquad \mathcal Q(v_n)\leq-\frac\kappa\eps,
	\qquad \delta_n\longrightarrow0.
	\]
	Next, we consider the set
	\[
	\widetilde{\mathcal N}_{\eps,\kappa}^-
	:=\left\{(t,w)\in[\frac{1}{2},2]\times\cS_{\eps\rho}\;\middle|\;
	\widetilde P(t,w)=0,\quad
	\widetilde{\mathcal Q}(t,w)\leq-\frac\kappa{2\eps}\right\}.
	\]
	Equip this set with the metric
	\[
	d_\times((t,w),(s,v))=|t-s|+\|w-v\|_{\cH^2}.
	\]
	The set $\widetilde{\mathcal N}_{\eps,\kappa}^-$ is closed in the complete
	metric space $[\frac{1}{2},2]\times\cS_{\eps\rho}$. Therefore,
	$(\widetilde{\mathcal N}_{\eps,\kappa}^-,d_\times)$ is a complete metric space.
	For each $(t,w)$ in this set,
	$t\star w\in\mathcal P_{\eps\rho}^-$, so we have $\widetilde E(t,w)\geq d_\eps$. Conversely, $(1,v_n)$ belongs to the
	set for $n$ large enough, with $\widetilde E(1,v_n)=E(v_n)\to d_\eps$.
	Hence $\inf\widetilde E=d_\eps$ on this complete metric space.
	
	Now Ekeland's variational principle \cite{Ekeland1974}, applied at $(1,v_n)$,
	gives $(t_n,w_n)$ in this set such that
	\begin{align*}
		& |t_n-1|+\|w_n-v_n\|_{\cH^2}\leq\sqrt{\delta_n},\\
		& \widetilde E(t_n,w_n)\leq E(v_n),\\
		& \widetilde E(t,w)\geq\widetilde E(t_n,w_n)
		-\sqrt{\delta_n}\bigl(|t-t_n|+\|w-w_n\|_{\cH^2}\bigr)
		\quad\text{for all }(t,w)\in\widetilde{\mathcal N}_{\eps,\kappa}^-.
	\end{align*}
	Set $u_n=t_n\star w_n$. Since scaling preserves signs and the $L^2$ norm,
	we obtain
	\[
	\sum_{i=1}^2\|(u_{i,n})_-\|_2^2
	=\sum_{i=1}^2\|(w_{i,n})_-\|_2^2
	\leq\|w_n-v_n\|_2^2\leq\delta_n\to0.
	\]
	Thus \eqref{eq:high-PS-asymptotic-nonnegative} holds. Since $t_n\to1$,
	we have $t_n\in(\frac{1}{2},2)$ for all sufficiently large $n$.
	Moreover, the energies satisfy
	\[
	d_\eps\leq E(u_n)=\widetilde E(t_n,w_n)\leq E(v_n)\longrightarrow d_\eps.
	\]
	Therefore, applying \cref{lem:fold-gap} to $u_n$ gives
	\[
	\begin{aligned}
	& \widetilde{\mathcal Q}(t_n,w_n)=\mathcal Q(u_n)
	\leq-\frac\kappa\eps<-\frac{\kappa}{2\eps},\\
	& \partial_t\widetilde P(t_n,w_n)
	=\frac{\widetilde{\mathcal Q}(t_n,w_n)}{t_n}\neq0.
	\end{aligned}
	\]
	The three equality constraints are $\widetilde P(t,w)=0$ and
	$M_i(w)=\frac12\|w_i\|_2^2=\frac{\eps\rho_i}{2}$ for $i=1,2$.
	Here $M_i$ depends only on $w$, and $DM_1(w_n)$ and $DM_2(w_n)$
	are linearly independent. Since $\partial_t\widetilde P(t_n,w_n)\neq0$,
	the derivatives of all three equality constraints are linearly
	independent. Define
	\[
	\begin{aligned}
		& X_n:=\{\xi=(s,h)\in\R\times\cH^2\mid D\widetilde P(t_n,w_n)[\xi]=0,\\
		& DM_i(w_n)[h]=0,\ i=1,2\}.
	\end{aligned}
	\]
	Let $\xi=(s,h)\in X_n$. By the implicit function theorem, there exists
	a continuously differentiable map $r\mapsto\gamma(r)$ with
	$\gamma(0)=(t_n,w_n)$ and $\gamma'(0)=\xi$.
	By construction, $\gamma(r)$ satisfies the three equality constraints
	for small $|r|$.
	Since $\frac{1}{2}<t_n<2$ and
	$\widetilde{\mathcal Q}(t_n,w_n)<-\frac{\kappa}{2\eps}$, continuity gives
	$\gamma(r)\in\widetilde{\mathcal N}_{\eps,\kappa}^-$ for sufficiently small $|r|$.
	For $r>0$, we apply Ekeland's inequality at $\gamma(r)$.
	After subtracting $\widetilde E(t_n,w_n)$ and dividing by $r$, we let
	$r\to0^+$. Using $\gamma'(0)=\xi=(s,h)$, we obtain
	\[
	D\widetilde E(t_n,w_n)[\xi]
	\geq-\sqrt{\delta_n}\bigl(|s|+\|h\|_{\cH^2}\bigr).
	\]
	Next, we apply the same argument to $r\mapsto\gamma(-r)$, whose
	derivative at zero is $-\xi$. Combining the two bounds, we obtain
	\[
	\bigl|D\widetilde E(t_n,w_n)[\xi]\bigr|
	\le
	\sqrt{\delta_n}\bigl(|s|+\|h\|_{\cH^2}\bigr)
	\le
	\sqrt{2\delta_n}\,\|\xi\|_{\mathbb R\times \cH^2}.
	\]
	Taking the supremum over $\xi\in X_n$ with
	$\|\xi\|_{\R\times\cH^2}=1$, we obtain
	\[
	\left\|D\widetilde E(t_n,w_n)\big|_{X_n}\right\|
	\leq\sqrt{2\delta_n}.
	\]
	The Hahn-Banach theorem gives a functional
	$R_n\in(\R\times\cH^2)'$ such that
	\[
	R_n[\xi]=D\widetilde E(t_n,w_n)[\xi]\quad\text{for all }\xi\in X_n,
	\qquad
	\|R_n\|_{(\R\times\cH^2)'}=
	\left\|D\widetilde E(t_n,w_n)\big|_{X_n}\right\|
	\leq\sqrt{2\delta_n}.
	\]
	Thus $D\widetilde E(t_n,w_n)-R_n$ vanishes on $X_n$.
	Since $X_n$ is the common kernel of the three linearly independent
	constraint differentials, every continuous linear functional that
	vanishes on $X_n$ is a linear combination of these differentials.
	Hence there exist multipliers
	$\eta_n,\lambda_{1,n},\lambda_{2,n}\in\R$ such that
	\begin{equation}\label{eq:augmented-lagrange}
		\begin{aligned}
			& 
			D\widetilde E(t_n,w_n)-\eta_nD\widetilde P(t_n,w_n) -\sum_{i=1}^2\lambda_{i,n}(0,DM_i(w_n))=R_n,\\
			& \|R_n\|_{(\R\times\cH^2)'}\leq\sqrt{2\delta_n}.
		\end{aligned}
	\end{equation}
	Next, we test \eqref{eq:augmented-lagrange} with $(1,0)$. Using
	$\partial_t\widetilde E(t_n,w_n)=\frac{\widetilde P(t_n,w_n)}{t_n}=0$,
	we obtain
	\[
	-\eta_n\frac{\widetilde{\mathcal Q}(t_n,w_n)}{t_n}
	=R_n[(1,0)],\qquad |R_n[(1,0)]|\leq\sqrt{2\delta_n}\to0.
	\]
	Since $\eps$ is fixed and
	$|\widetilde{\mathcal Q}(t_n,w_n)|\geq\kappa/\eps$,
	it follows that $\eta_n\to0$.
	
	These estimates prove \eqref{eq:high-PS-level}. Moreover, the identity
	$K(u_n)+5W(u_n)=6E(u_n)$, together with $E(u_n)\to d_\eps$ and the
	fixed masses, shows that $\{u_n\}$ is bounded in $\cH^2$.
	To express the derivative estimate in the original variables, take $h\in\cH^2$ and set $k=t_n^{-1}\star h$. Linearity of scaling
	and preservation of the $L^2$ inner product give
	\begin{align*}
		& D_w\widetilde E(t_n,w_n)[k] =E'(u_n)[h],\\
		& D_w\widetilde P(t_n,w_n)[k] =P'(u_n)[h],\\
		& DM_i(w_n)[k] =DM_i(u_n)[h].
	\end{align*}
	Since $t_n\to1$, the scaling maps and their inverses are uniformly
	bounded on $\cH^2$. Moreover, the sequence $P'(u_n)$ is bounded in $(\cH^2)'$.
	Now we apply \eqref{eq:augmented-lagrange} to $(0,k)$.
	Together with the preceding bounds and $\eta_n\to0$, the resulting
	identity gives \eqref{eq:true-product-PS}.
	
	Testing \eqref{eq:true-product-PS} with the bounded pair $u_{i,n}e_i$
	gives, for $i=1,2$,
	\[
	\lambda_{i,n}\eps\rho_i
	=\int_{\R^3}(|\nabla u_{i,n}|^2+|y|^2u_{i,n}^2)\dd x
	-\mu_iN_i(u_n)-2\beta C(u_n)+o(1).
	\]
	The $\cH^2$ bound and the Sobolev embedding control the quadratic and
	quartic terms. Since $\eps\rho_i>0$ is fixed, both multiplier
	sequences are bounded.
\end{proof}

Before extracting profiles from this sequence, we show that each
nonzero component of a nonnegative solution of \eqref{eq:system} is
strictly positive and has a multiplier satisfying $\lambda_i<\Lambda_0$.
To estimate the Lagrange multipliers, we test against the transverse
ground state as in
\cite[Proposition~2.6]{JeanjeanSong2025}.

\begin{lemma}\label{lem:profile-equations-pohozaev}
	Let $U\in\cH^2$ be a nonzero, componentwise nonnegative weak solution
	of \eqref{eq:system} for some $\lambda_1,\lambda_2\in\R$.
	Every nonzero component $U_i$ is strictly positive in $\R^3$.
	Moreover,
	\begin{equation}\label{eq:profile-P-zero}
		P(U)=0,\qquad E(U)=\frac{K(U)+5W(U)}6>0.
	\end{equation}
	For each $i$ with $U_i\not\equiv0$, the corresponding Lagrange
	multiplier satisfies
	\begin{equation}\label{eq:profile-subthreshold-multiplier}
		\lambda_i<\Lambda_0=2.
	\end{equation}
\end{lemma}

\begin{proof}
	Let $\{i,k\}=\{1,2\}$. The embedding $H^1(\R^3)\hookrightarrow L^6(\R^3)$
	and H\"older's inequality give
	\[
	\|U_i^3\|_2=\|U_i\|_6^3,\qquad
	\|U_iU_k^2\|_2\leq\|U_i\|_6\|U_k\|_6^2.
	\]
	Since the nonlinearities belong to $L^2(\R^3)$, the equation has a locally
	$L^2$ right-hand side. Local elliptic regularity yields
	$U_i\in W^{2,2}_{\mathrm{loc}}(\R^3)$. Hence
	$U_i\in C^{0,\alpha}_{\mathrm{loc}}(\R^3)$ for every
	$\alpha\in(0,\tfrac12)$. The right-hand side is then locally H\"older
	continuous, and Schauder regularity gives
	$U_i\in C^{2,\alpha}_{\mathrm{loc}}(\R^3)$.
	Fix $R>0$ and choose
	$$c_R\geq\max\{0,\sup_{B_R}(|y|^2-\lambda_i)\}.$$
	 Then
	\[
	-\Delta U_i+c_RU_i
	=(c_R+\lambda_i-|y|^2)U_i+\mu_iU_i^3+2\beta U_iU_k^2\geq0
	\quad\text{in }B_R.
	\]
	The strong maximum principle shows that $U_i$ is either identically zero
	or strictly positive throughout $\R^3$.
	
	Testing the $i$th equation in \eqref{eq:system} with $U_i$ gives
	\[
	\int_{\R^3}\bigl(|\nabla U_i|^2+|y|^2U_i^2\bigr)\dd x
	=\lambda_i\|U_i\|_2^2+\mu_iN_i(U)+2\beta C(U),
	\qquad i=1,2.
	\]
	Summing these identities and recalling
	$4B(U)=\mu_1N_1(U)+\mu_2N_2(U)+4\beta C(U)$ yields
	\begin{equation}\label{eq:profile-nehari}
		K(U)+W(U)=\sum_{i=1}^2\lambda_i\|U_i\|_2^2+4B(U).
	\end{equation}
	If both components are nonzero, the Pohozaev identity in
	\cite[Lemma~2.9]{JiaLiLuo2020}, with coupling coefficient $2\beta$,
	gives $K(U)-W(U)=3B(U)$.
	Suppose $U_k\equiv0$ for some $k\in\{1,2\}$. Set $i:=3-k$ and
	$v:=U_i$. Then $v$ satisfies the scalar equation
	\[
	-\Delta v+|y|^2v=\lambda_i v+\mu_i v^3\quad\text{in }\R^3.
	\]
	Set
	\[
	K_v=\int_{\R^3}|\nabla v|^2\dd x,\quad
	W_v=\int_{\R^3}|y|^2v^2\dd x,\quad
	M_v=\|v\|_2^2,\quad N_v=\|v\|_4^4.
	\]
	Testing the scalar equation with $v$ and applying the spatial Pohozaev identity give
	\begin{align}
		& K_v+W_v =\lambda_i M_v+\mu_i N_v,
		\label{eq:scalar-profile-nehari}\\
		& \frac12K_v+\frac52W_v =\frac32\lambda_i M_v+\frac34\mu_i N_v.
		\label{eq:scalar-profile-pohozaev}
	\end{align}
	To prove the second identity, choose a smooth function $\chi$ that
	equals one on $B_1$ and vanishes outside $B_2$.
	Set $\chi_R(x)=\chi(x/R)$ and test the scalar equation with
	$\chi_R x\cdot\nabla v$.
	The local regularity of $v$ justifies integration by parts.
	Set $A_R:=B_{2R}\setminus B_R$. Since
	$\operatorname{supp}\nabla\chi_R\subset A_R$ and
	$|x||\nabla\chi_R|\leq C$, the total absolute contribution from the
	region where $\chi_R$ varies is bounded by
	\[
	C\int_{A_R}\left(|\nabla v|^2+(1+|y|^2)v^2+|v|^4\right)\dd x
	\longrightarrow0\quad\text{as }R\to\infty,
	\]
	since $v\in\cH\cap L^4(\R^3)$.
	Using $x\cdot\nabla(|y|^2)=2|y|^2$ in the integration by parts yields
	\eqref{eq:scalar-profile-pohozaev}.
	Subtracting $\frac32$ times \eqref{eq:scalar-profile-nehari} from
	\eqref{eq:scalar-profile-pohozaev} yields
	$K_v-W_v=\frac34\mu_iN_v=3B(U)$.
	Thus $P(U)=0$ holds for semitrivial profiles. Now the energy identity in
	\eqref{eq:profile-P-zero} follows from \eqref{eq:sheet-identities}.
	The strict positivity of the energy follows from \eqref{eq:spectral-bound} and $U\neq0$.
	
	Now fix $i\in\{1,2\}$ such that $U_i\not\equiv0$ and define
	\begin{align*}
		& h_i(z) =\int_{\R^2}\phi_0(y)U_i(y,z)\dd y,\\
		& G_i(z) =\int_{\R^2}\phi_0(y)
		(\mu_iU_i^3+2\beta U_iU_k^2)(y,z)\dd y.
	\end{align*}
	The Cauchy-Schwarz inequality gives $h_i\in H^1(\R)$ and $G_i\in L^2(\R)$.
	Given $\psi\in C_c^\infty(\R)$, we test the $i$th equation in
	\eqref{eq:system} with $\phi_0(y)\psi(z)\in\cH$.
	Using \eqref{eq:transverse-ground-state}, we obtain
	\[
	\int_\R h_i'\psi'\dd z
	+(\Lambda_0-\lambda_i)\int_\R h_i\psi\dd z
	=\int_\R G_i\psi\dd z.
	\]
	Equivalently,
	\begin{equation}\label{eq:ground-channel-profile}
		-h_i''+(\Lambda_0-\lambda_i)h_i=G_i
		\quad\text{in }\mathcal D'(\R).
	\end{equation}
	Here $h_i\geq0$ is nonzero because $\phi_0>0$ and $U_i\geq0$ is nonzero.
	Also, $G_i\geq0$. Equation \eqref{eq:ground-channel-profile} gives
	$h_i\in H^2(\R)$. If $\lambda_i\geq\Lambda_0$, then $h_i''\leq0$.
	Concavity implies
	$h_i(t)\leq h_i(s)+h_i'(s)(t-s)$ for all $s,t\in\R$.
	Since $h_i\geq0$ on $\R$, necessarily $h_i'(s)=0$ for every $s$.
	Thus $h_i$ is constant. Since $h_i\in L^2(\R)$, the constant must be zero, a contradiction.
	This proves \eqref{eq:profile-subthreshold-multiplier}.
\end{proof}

Now we establish an axial profile decomposition.
Fix $\eps\in(0,\eps_0)$ and let $\{u_n\}$ be the sequence constructed in
\cref{prop:high-PS}.
The Lagrange multiplier sequences $\{\lambda_{1,n}\}$ and
$\{\lambda_{2,n}\}$ are bounded.
Consequently, up to a subsequence, there exist $\lambda_1,\lambda_2\in\R$
such that
\begin{equation}\label{eq:high-lambda-limit}
	\lambda_{i,n}\to\lambda_i\in\R\quad(i=1,2),\qquad n\to\infty.
\end{equation}
The common axial translations $\tau_s$ introduced in
\cref{sec:introduction} preserve the $\cH^2$ norm and the values of
$K,W,B$. For every fixed $u\in\cH^2$,
\[
\tau_su\rightharpoonup0
\quad\text{in }\cH^2\quad\text{as }|s|\to\infty.
\]
To prove the weak convergence, first take
$u,v\in C_c^\infty(\R^3;\R^2)$.
Since the supports of $\tau_su$ and $v$ are disjoint for sufficiently
large $|s|$, it follows that
\[
\langle\tau_su,v\rangle_{\cH^2}=0.
\]
By density and preservation of the $\cH^2$ norm under $\tau_s$, we obtain
\[
\langle\tau_su,v\rangle_{\cH^2}\longrightarrow0
\quad\text{as }|s|\to\infty
\]
for all $u,v\in\cH^2$.

\begin{lemma}\label{lem:axial-profile-decomposition}
	Let $\{u_n\}$ be the sequence in \cref{prop:high-PS}.
	Up to a subsequence, there exist an integer $J\geq1$, nonzero profiles
	$U^1,\ldots,U^J\in\cH^2$, and sequences $\{z_n^j\}\subset\R$
	such that
	\[
	|z_n^j-z_n^k|\to\infty\quad(1\leq j\neq k\leq J),\qquad
	r_n=u_n-\sum_{j=1}^J\tau_{z_n^j}U^j.
	\]
	Each profile $U^j=(U_1^j,U_2^j)$ is a nonnegative weak solution
	of the following system in $\R^3$:
	\begin{equation}\label{eq:profile-limit-system}
		-\Delta U_i^j+|y|^2U_i^j
		=\lambda_iU_i^j+\mu_i(U_i^j)^3+2\beta U_i^j(U_k^j)^2,
		\qquad \{i,k\}=\{1,2\}.
	\end{equation}
	Here $\lambda_1,\lambda_2$ are the limits in
	\eqref{eq:high-lambda-limit} and are the same for all profiles.
	Moreover, $P(U^j)=0$ and $E(U^j)>0$ for every $j$.
	If $U_i^j\not\equiv0$, then $U_i^j>0$ in $\R^3$ and
	$\lambda_i<\Lambda_0$.
	In addition, the following decompositions hold as $n\to\infty$:
	\begin{align}
		& \|u_{i,n}\|_2^2 =\sum_{j=1}^J\|U_i^j\|_2^2+\|r_{i,n}\|_2^2+o(1),
		\qquad i=1,2,\label{eq:profile-mass-split}\\
		& K(u_n)+W(u_n) =\sum_{j=1}^J\bigl(K(U^j)+W(U^j)\bigr)
		+K(r_n)+W(r_n)+o(1),
		\label{eq:profile-quadratic-split}\\
		& B(u_n) =\sum_{j=1}^J B(U^j)+B(r_n)+o(1).
		\label{eq:profile-B-split}
	\end{align}
	Moreover, the remainder satisfies
	\begin{equation}\label{eq:profile-remainder-L4}
		\|r_n\|_4\longrightarrow0.
	\end{equation}
	Finally, for $j=1,\ldots,J$, the translated sequence satisfies
	\[
	\begin{aligned}
		& \tau_{-z_n^j}u_n \rightharpoonup U^j & & \text{in }\cH^2,\\
		& \tau_{-z_n^j}u_n \longrightarrow U^j & & \text{in }L^p_{\mathrm{loc}},
		\qquad 2\leq p<6.
	\end{aligned}
	\]
\end{lemma}

\begin{proof}
	\noindent\textbf{Step 1.}
	First, we show that a bounded sequence in $\cH^2$ that does not
	converge to zero in $L^4$ admits a nonzero axial profile up to a subsequence.
	Let $Q_k=k+(0,1)^3$, $k\in\mathbb Z^3$.
	H\"older's inequality and the Sobolev embedding on unit cubes give
	$\|f\|_{L^4(Q_k)}^4\leq C\|f\|_{L^2(Q_k)}\|f\|_{H^1(Q_k)}^3$,
	with $C$ independent of $k$.
	Summing this estimate and using
	$\sum_k a_k^3\leq(\sum_k a_k^2)^{\frac{3}{2}}$ for $a_k\geq0$
	gives, for some $C_{\rm loc}>0$,
	\begin{equation}\label{eq:local-concentration-estimate}
		\|f\|_4^4\leq C_{\rm loc}
		\left(\sup_{k\in\mathbb Z^3}\|f\|_{L^2(Q_k)}\right)
		\|f\|_{H^1(\R^3)}^3
		\qquad(f\in H^1(\R^3)).
	\end{equation}
	For $v=(v_1,v_2)\in H$, the chain rule for weak derivatives gives
	\[
	\bigl\||v|\bigr\|_{H^1(\R^3)}^2
	\leq\sum_{i=1}^2\|v_i\|_{H^1(\R^3)}^2.
	\]
	Therefore, applying \eqref{eq:local-concentration-estimate} to $|v|$
	gives the corresponding estimate for pairs.
	This local concentration argument is used in \cite{Lions1984I}.
	
	Let $\{v_n\}$ be bounded in $\cH^2$ and suppose that
	$\limsup_{n\to\infty}\|v_n\|_4>0$.
	Up to a subsequence, \eqref{eq:local-concentration-estimate} gives
	$p_n\in\mathbb Z^2$, $q_n\in\mathbb Z$, and $\eta>0$ such that
	\[
	\int_{Q_{(p_n,q_n)}}|v_n|^2\dd x\geq\eta.
	\]
	If $|p_n|>\sqrt2$, then every $(y,z)\in Q_{(p_n,q_n)}$
	satisfies $|y|\geq|p_n|-\sqrt2$. It follows that
	\[
	\int_{Q_{(p_n,q_n)}}|v_n|^2\dd x
	\leq\frac{W(v_n)}{(|p_n|-\sqrt2)^2}.
	\]
	Since $W(v_n)$ is bounded and $\eta>0$, the sequence $\{p_n\}$
	takes only finitely many values. Consequently, up to a subsequence,
	we may assume $p_n=p$ for some $p\in\mathbb Z^2$.
	The sequence $\{\tau_{-q_n}v_n\}$ is bounded in $\cH^2$.
	By weak compactness and the Rellich-Kondrachov theorem
	\cite[Theorem~9.16]{Brezis2011}, we may assume, up to a subsequence, that
	\[
	\begin{aligned}
		& \tau_{-q_n}v_n\rightharpoonup V &&\text{in }\cH^2,\\
		& \tau_{-q_n}v_n\longrightarrow V &&\text{in }L^p_{\mathrm{loc}},\quad 2\leq p<6,\\
		& \tau_{-q_n}v_n\longrightarrow V &&\text{almost everywhere in }\R^3.
	\end{aligned}
	\]
	In particular, convergence in $L^2(Q_{(p,0)};\R^2)$ gives
	\[
	\int_{Q_{(p,0)}}|V|^2\dd x
	=\lim_{n\to\infty}\int_{Q_{(p,0)}}|\tau_{-q_n}v_n|^2\dd x
	=\lim_{n\to\infty}\int_{Q_{(p,q_n)}}|v_n|^2\dd x\geq\eta.
	\]
	Therefore, $V\neq0$.
	
	We apply this argument to $\{u_n\}$.
	Since $P(u_n)=0$, \eqref{eq:high-PS-level} gives
	\[
	3B(u_n)=4W(u_n)-\mathcal Q(u_n)\geq\frac\kappa{2\eps}.
	\]
	Then the bound $B(u_n)\leq C\|u_n\|_4^4$ implies
	$\liminf_n\|u_n\|_4>0$, so the preceding argument yields a nonzero profile.
	
	\smallskip\noindent\textbf{Step 2.}
	Next, we construct successive profiles and prove the separation of
	their translation parameters.
	Set $r_n^0=u_n$.
	Since $\{u_n\}$ is bounded and axial translations preserve the norm,
	each fixed remainder sequence $\{r_n^{j-1}\}$ is bounded in $\cH^2$.
	If $\limsup_{n\to\infty}\|r_n^{j-1}\|_4>0$, we apply Step~1 to
	$\{r_n^{j-1}\}$. Up to a subsequence, there exist $U^j\neq0$ and
	$\{z_n^j\}\subset\R$ such that
	\[
	\tau_{-z_n^j}r_n^{j-1}\rightharpoonup U^j\quad\text{in }\cH^2.
	\]
	Moreover, the convergence holds in $L^p_{\mathrm{loc}}$ for $2\leq p<6$
	and almost everywhere.
	Setting $r_n^j=r_n^{j-1}-\tau_{z_n^j}U^j$, we obtain
	$\tau_{-z_n^j}r_n^j\rightharpoonup0$ in $\cH^2$.
	Next, we verify the separation property
	\[
	|z_n^j-z_n^k|\to\infty\qquad(1\leq k<j).
	\]
	For $j>1$, assume that
	$\tau_{-z_n^k}r_n^{j-1}\rightharpoonup0$ for $1\leq k<j$.
	If $|z_n^j-z_n^k|$ did not tend to infinity for some $k<j$,
	then up to a subsequence, $z_n^j-z_n^k\to s\in\R$.
	For every fixed $\varphi\in\cH^2$, strong continuity of translations gives
	$\tau_{z_n^j-z_n^k}\varphi\to\tau_s\varphi$ in $\cH^2$.
	Together with the preceding weak convergence, this yields
	\[
	\bigl\langle\tau_{-z_n^j}r_n^{j-1},\varphi\bigr\rangle_{\cH^2}
	=\bigl\langle\tau_{-z_n^k}r_n^{j-1},
	\tau_{z_n^j-z_n^k}\varphi\bigr\rangle_{\cH^2}\longrightarrow0,
	\]
	contrary to $U^j\neq0$.
	For $1\leq k<j$, the separation just proved gives
	\[
	\tau_{-z_n^k}r_n^j
	=\tau_{-z_n^k}r_n^{j-1}-\tau_{z_n^j-z_n^k}U^j
	\rightharpoonup0.
	\]
	Together with the case $k=j$, this establishes the same weak convergence
	for $k=1,\ldots,j$, so the argument can be repeated.
	
	\smallskip\noindent\textbf{Step 3.}
	Now we show that each profile is a translated weak limit of the
	original sequence and satisfies the limiting system. Indeed,
	we have
	\[
	\tau_{-z_n^j}u_n
	=\tau_{-z_n^j}r_n^{j-1}
	+\sum_{k=1}^{j-1}\tau_{z_n^k-z_n^j}U^k.
	\]
	By the separation property, every term in the sum converges weakly to zero.
	Consequently, $\tau_{-z_n^j}u_n\rightharpoonup U^j$ in $\cH^2$.
	Moreover, by the Rellich-Kondrachov theorem, we may assume,
	up to a subsequence, that
	\begin{equation}\label{eq:profile-local-convergence}
		\tau_{-z_n^j}u_n\longrightarrow U^j
		\quad\text{in }L^p_{\mathrm{loc}},
		\qquad 2\leq p<6.
	\end{equation}
	The local $L^2$ convergence and
	\eqref{eq:high-PS-asymptotic-nonnegative} imply
	$U_i^j\geq0$ for $i=1,2$.
	
	To derive the equation satisfied by $U^j$, we use the derivative
	estimate for the original Palais-Smale sequence.
	By \eqref{eq:true-product-PS}, the residual $F_n\in(\cH^2)'$
	satisfies
	\begin{equation}\label{eq:profile-PS-residual}
		\begin{aligned}
			& F_n :=E'(u_n)-\sum_{i=1}^2\lambda_{i,n}u_{i,n}e_i,\\
			& \|F_n\|_{(\cH^2)'} \longrightarrow0.
		\end{aligned}
	\end{equation}
	Fix $\varphi=(\varphi_1,\varphi_2)$ with
	$\varphi_1,\varphi_2\in C_c^\infty(\R^3)$ and apply $F_n$ to
	$\tau_{z_n^j}\varphi$. Since axial translations preserve the
	$\cH^2$ norm, \eqref{eq:profile-PS-residual} gives
	\[
	|F_n[\tau_{z_n^j}\varphi]|
	\leq\|F_n\|_{(\cH^2)'}\|\varphi\|_{\cH^2}\longrightarrow0.
	\]
	After the change of variables $z=\zeta+z_n^j$, weak convergence
	of $\tau_{-z_n^j}u_n$ in $\cH^2$ and
	$\lambda_{i,n}\to\lambda_i$ allow us to pass to the limit in
	the gradient, potential, and Lagrange multiplier terms.
	By \eqref{eq:profile-local-convergence}, we have
	\[
	\tau_{-z_n^j}u_n\longrightarrow U^j
	\quad\text{in }L^4_{\mathrm{loc}}.
	\]
	Together with H\"older's inequality, this convergence yields, for
	$\{i,k\}=\{1,2\}$,
	\[
	\begin{aligned}
		& (\tau_{-z_n^j}u_{i,n})^3 \longrightarrow(U_i^j)^3,\qquad\text{in }L_{\mathrm{loc}}^{\frac{4}{3}},\\
		& (\tau_{-z_n^j}u_{i,n})(\tau_{-z_n^j}u_{k,n})^2 \longrightarrow U_i^j(U_k^j)^2,\qquad\text{in }L_{\mathrm{loc}}^{\frac{4}{3}}.
	\end{aligned}
	\]
	These limits and \eqref{eq:profile-PS-residual} establish
	the weak formulation of \eqref{eq:profile-limit-system}
	for smooth compactly supported test functions.
	By density and continuity of the terms, this identity holds for
	every test function in $\cH$.
	Hence \cref{lem:profile-equations-pohozaev} gives the positivity,
	Pohozaev identity, and multiplier estimates stated in the lemma.
	
	\smallskip\noindent\textbf{Step 4.}
	Next, we prove the mass and functional decompositions for any fixed
	number of profiles.
	At each step, expanding the squared norms and using
	$\tau_{-z_n^j}r_n^j\rightharpoonup0$ gives
	\[
	\begin{aligned}
		& \|r_{i,n}^{j-1}\|_2^2 =\|U_i^j\|_2^2+\|r_{i,n}^j\|_2^2+o(1),\qquad i=1,2,\\
		& K(r_n^{j-1})+W(r_n^{j-1}) =K(U^j)+W(U^j)+K(r_n^j)+W(r_n^j)+o(1).
	\end{aligned}
	\]
	The scalar Br\'ezis-Lieb lemma \cite{BrezisLieb1983} gives the
	decompositions of $N_1$ and $N_2$.
	By Step~2 and the Sobolev embedding,
	\[
	\begin{aligned}
	&\tau_{-z_n^j}r_n^{j-1}\longrightarrow U^j
	\quad\text{strongly in }L^4_{\mathrm{loc}}(\R^3),\\
	&\tau_{-z_n^j}r_n^{j-1}\longrightarrow U^j
	\quad\text{almost everywhere in }\R^3,\\
	&\sup_n\|\tau_{-z_n^j}r_n^{j-1}\|_4<\infty.
	\end{aligned}
	\]
	The mixed Br\'ezis-Lieb lemma \cite[Lemma~2.4]{GouJeanjean2018},
	with $r_1=r_2=2$, therefore yields
	\[
	C(r_n^{j-1})=C(U^j)+C(r_n^j)+o(1).
	\]
	Hence $B(r_n^{j-1})=B(U^j)+B(r_n^j)+o(1)$.
	At each step, these identities hold up to a subsequence
	for which all earlier limits remain valid.
	Summing the identities after $N$ steps gives
	\eqref{eq:profile-mass-split}--\eqref{eq:profile-B-split}
	with $J=N$ and $r_n=r_n^N$.
	In particular,
	\begin{equation}\label{eq:finite-profile-energy-sum}
		\sum_{j=1}^N\bigl(K(U^j)+W(U^j)\bigr)
		\leq\limsup_{n\to\infty}\bigl(K(u_n)+W(u_n)\bigr).
	\end{equation}
	
	\smallskip\noindent\textbf{Step 5.}
	Finally, we prove that only finitely many nonzero profiles can occur.
	For this purpose, we establish a positive lower bound for
	$K(U^j)+W(U^j)$ independent of $j$.
	Let $I_\lambda=\{i\in\{1,2\}\mid\lambda_i<\Lambda_0\}$ and set
	\[
	\lambda_+=\max\bigl(\{0\}\cup\{\lambda_i\mid i\in I_\lambda\}\bigr),
	\qquad\theta_\lambda=1-\frac{\lambda_+}{\Lambda_0}.
	\]
	By \cref{lem:profile-equations-pohozaev},
	$I_\lambda\neq\varnothing$ and $U_i^j=0$ for $i\notin I_\lambda$.
	Thus $\theta_\lambda>0$ is independent of $j$.
	The spectral estimate \eqref{eq:spectral-bound} gives
	\[
	\sum_{i=1}^2\|U_i^j\|_2^2
	\leq\frac{K(U^j)+W(U^j)}{\Lambda_0}.
	\]
	It follows that
	\[
	\sum_{i=1}^2\lambda_i\|U_i^j\|_2^2
	\leq\frac{\lambda_+}{\Lambda_0}\bigl(K(U^j)+W(U^j)\bigr).
	\]
	Substitution into \eqref{eq:profile-nehari} yields
	\begin{equation}\label{eq:profile-lower-nonlinearity}
			 \theta_\lambda\bigl(K(U^j)+W(U^j)\bigr) \leq K(U^j)+W(U^j)-\sum_{i=1}^2\lambda_i\|U_i^j\|_2^2=4B(U^j).
	\end{equation}
	Moreover, the spectral estimate gives
	$\sum_i\|U_i^j\|_{H^1}^2\leq(1+\Lambda_0^{-1})(K(U^j)+W(U^j))$.
	Then the Sobolev embedding and H\"older's inequality provide a constant
	$C_4>0$, independent of $j$, such that
	\[
	4B(U^j)\leq C_4\bigl(K(U^j)+W(U^j)\bigr)^2.
	\]
	Combining this bound with \eqref{eq:profile-lower-nonlinearity}
	and using $U^j\neq0$ gives
	\begin{equation}\label{eq:profile-quadratic-quantum}
		K(U^j)+W(U^j)\geq\frac{\theta_\lambda}{C_4}>0.
	\end{equation}
	Since $P(u_n)=0$, \eqref{eq:sheet-identities} gives
	$K(u_n)+W(u_n)\leq6E(u_n)$.
	Together with $E(u_n)\to d_\eps$ and
	\eqref{eq:finite-profile-energy-sum}, this yields
	\[
	N\frac{\theta_\lambda}{C_4}
	\leq\sum_{j=1}^N\bigl(K(U^j)+W(U^j)\bigr)\leq6d_\eps.
	\]
	Thus $N\leq\frac{6C_4d_\eps}{\theta_\lambda}$ for every finite number $N$
	of profiles obtained above. At each finite stage, if
	$\limsup_{n\to\infty}\|r_n^N\|_4>0$, Step~2 applies to a further
	subsequence, preserves the profiles already extracted, and produces
	another nonzero profile; the uniform bound on the number of profiles
	therefore yields a finite $J$ such that $\|r_n^J\|_4\to0$.
	Step~1 ensures that $J\geq1$.
	Setting $r_n=r_n^J$ gives
	\eqref{eq:profile-mass-split}--\eqref{eq:profile-remainder-L4}
	and completes the proof.
\end{proof}

Next, we identify the limits of the component masses and the quadratic
energy of $r_n$. Define
\begin{align}
	& b_i =\eps\rho_i-\sum_{j=1}^J\|U_i^j\|_2^2\geq0, & & b =(b_1,b_2),\label{eq:missing-profile-mass}\\
	& \ell =d_\eps-\sum_{j=1}^JE(U^j),\label{eq:ell-profile}\\
	& L_0(b) =\frac{\Lambda_0}{2}(b_1+b_2).\label{eq:L0-definition}
\end{align}
The nonnegativity of $b_i$ follows from
\eqref{eq:profile-mass-split}.

\begin{lemma}\label{lem:residual-ledger}
	Let $\{r_n\}$ be the remainder sequence in
	\cref{lem:axial-profile-decomposition}. Then the following limits
	hold as $n\to\infty$:
	\[
	\begin{aligned}
		& \|r_{i,n}\|_2^2 \to b_i,\qquad i=1,2,\\
		& \frac12\bigl(K(r_n)+W(r_n)\bigr) \to\ell\geq L_0(b).
	\end{aligned}
	\]
	Moreover, the following convergence holds:
	\begin{equation}\label{eq:Q-ledger}
		\mathcal Q(u_n)\longrightarrow q:=\sum_{j=1}^J\mathcal Q(U^j)+4\ell
		\leq-\frac{\kappa}{2\eps}.
	\end{equation}
\end{lemma}

\begin{proof}
	Since $\|u_{i,n}\|_2^2=\eps\rho_i$,
	\eqref{eq:profile-mass-split} gives
	\[
	\|r_{i,n}\|_2^2
	=\eps\rho_i-\sum_{j=1}^J\|U_i^j\|_2^2+o(1)
	\longrightarrow b_i,
	\qquad i=1,2.
	\]
	By \eqref{eq:profile-remainder-L4}, $B(r_n)\to0$.
	The decompositions \eqref{eq:profile-quadratic-split} and
	\eqref{eq:profile-B-split}, together with the definition of $E$, give
	\[
	\begin{aligned}
		 \frac12\bigl(K(r_n)+W(r_n)\bigr) &=E(u_n)-\sum_{j=1}^JE(U^j)+B(r_n)+o(1)\\
		& \longrightarrow d_\eps-\sum_{j=1}^JE(U^j)=\ell.
	\end{aligned}
	\]
	Applying \eqref{eq:spectral-bound} to $r_n$ and passing to the limit using the mass
	convergence above gives
	\[
	\ell=\lim_{n\to\infty}\frac12\bigl(K(r_n)+W(r_n)\bigr)
	\geq\frac{\Lambda_0}{2}(b_1+b_2)=L_0(b).
	\]
	
	Finally, we combine
	\eqref{eq:profile-quadratic-split}--\eqref{eq:profile-B-split}
	with $\mathcal Q=2(K+W)-9B$ to obtain
	\[
	\begin{aligned}
		 \mathcal Q(u_n) &=\sum_{j=1}^J\mathcal Q(U^j)
		+2\bigl(K(r_n)+W(r_n)\bigr)-9B(r_n)+o(1)\\
		& \longrightarrow\sum_{j=1}^J\mathcal Q(U^j)+4\ell=q.
	\end{aligned}
	\]
	Then the bound in \eqref{eq:high-PS-level} gives
	$q\leq-\frac{\kappa}{2\eps}$, proving \eqref{eq:Q-ledger}.
\end{proof}

\section{Proof of Theorem 1.1}\label{sec:attainment}

In this section, we complete the proof of \cref{thm:main}.
The main task is to exclude splitting and loss of mass in the profile
decomposition obtained in \cref{sec:profiles}. To this end, we first
establish the following estimates for combining two profiles.

\begin{lemma}
	\label{lem:two-block-quadrature-separation}
	Let $f,g\in\cH^2$. For $s\in\R$, define
	\[
	h_i^s:=\bigl(f_i^2+(\tau_sg_i)^2\bigr)^{\frac{1}{2}},\qquad i\in\{1,2\}.
	\]
	Then $h^s\in\cH^2$ and
	\[
	\|h_i^s\|_2^2=\|f_i\|_2^2+\|g_i\|_2^2,
	\qquad W(h^s)=W(f)+W(g).
	\]
	Moreover,
	\begin{equation}\label{eq:two-block-energy-comparison}
		K(h^s)\leq K(f)+K(g),\qquad B(h^s)\geq B(f)+B(g),
	\end{equation}
	and
	\begin{equation}\label{eq:two-block-separation-limit}
		K(h^s)\to K(f)+K(g),\qquad B(h^s)\to B(f)+B(g)
		\quad\text{as }|s|\to\infty.
	\end{equation}
\end{lemma}

\begin{proof}
First, the definition of $h^s$ and a change of variables in $z$ give
the mass and potential identities. By the chain rule for weak derivatives,
we have
\[
|\nabla h_i^s|^2\leq|\nabla f_i|^2+|\nabla\tau_sg_i|^2.
\]
Together with the mass and potential identities, this shows that
$h^s\in\cH^2$ and $K(h^s)\leq K(f)+K(g)$.

Next, expanding the quartic terms gives
\[
\begin{aligned}
    &B(h^s)-B(f)-B(g)\\
    =& {} \sum_{i=1}^2\frac{\mu_i}{2}
    \int_{\R^3}f_i^2(\tau_sg_i)^2\dd x+\beta\int_{\R^3}
    \bigl(f_1^2(\tau_sg_2)^2+f_2^2(\tau_sg_1)^2\bigr)\dd x.
\end{aligned}
\]
Each integral on the right is nonnegative and tends to zero as
$|s|\to\infty$, since $f_i^2,g_k^2\in L^2(\R^3)$ and
$\tau_s(g_k^2)\rightharpoonup0$ in $L^2(\R^3)$ for $i,k=1,2$.
Consequently, the inequality for $B$ in \eqref{eq:two-block-energy-comparison}
and its limit in \eqref{eq:two-block-separation-limit} follow.

Finally, we prove the convergence of $K(h^s)$ as $|s|\to\infty$.
Fix $R>0$ and let $B_R=\{x\in\R^3\mid |x|<R\}$.
By the definition of $h^s$, we have
\[
\begin{aligned}
    & 0\leq h_i^s-|f_i|\leq|\tau_sg_i|,\\
    & 0\leq\tau_{-s}h_i^s-|g_i|\leq|\tau_{-s}f_i|,
    \qquad i=1,2.
\end{aligned}
\]
Since $f_i,g_i\in L^2(\R^3)$, the translated functions on the right
converge to zero in $L^2(B_R)$ as $|s|\to\infty$. It follows that
\[
h_i^s\longrightarrow|f_i|,\qquad
\tau_{-s}h_i^s\longrightarrow|g_i|
\quad\text{in }L^2(B_R),\qquad i=1,2.
\]
Moreover, the mass identities and the kinetic bound above imply
weak convergence to the same limits in $H^1(B_R)$.
For $|s|>2R$, the ball $B_R$ is disjoint from the ball of radius $R$
centered at $(0,0,s)$. Consequently, we obtain
\[
K(h^s)\geq\sum_{i=1}^2\left(
\int_{B_R}|\nabla h_i^s|^2\dd x+
\int_{B_R}|\nabla(\tau_{-s}h_i^s)|^2\dd x\right).
\]
By weak lower semicontinuity and the identities
$|\nabla|f_i||=|\nabla f_i|$ and
$|\nabla|g_i||=|\nabla g_i|$ almost everywhere, we obtain
\[
\liminf_{|s|\to\infty}K(h^s)
\geq\sum_{i=1}^2\int_{B_R}
\bigl(|\nabla f_i|^2+|\nabla g_i|^2\bigr)\dd x.
\]
Letting $R\to\infty$ and using $K(h^s)\leq K(f)+K(g)$, we conclude that
$K(h^s)\to K(f)+K(g)$ as $|s|\to\infty$.
This proves the remaining limit in \eqref{eq:two-block-separation-limit}.
\end{proof}

Next, we combine the profiles with a function carrying the missing
masses $b_i$. For this construction, define
\[
E_*:=\sum_{j=1}^J E(U^j)+L_0(b).
\]
By \cref{lem:residual-ledger} and \eqref{eq:ell-profile}, we obtain
\[
E_*\leq\sum_{j=1}^J E(U^j)+\ell=d_\eps.
\]

\begin{lemma}
	\label{lem:finite-quadrature-gluing}
	Let $U^j,b,q$ be as in
	\cref{lem:axial-profile-decomposition,lem:residual-ledger}.
	Assume $b\neq0$ or $J\geq2$. Then, for every $\delta>0$, there exist
	$w\in\cS_{\eps\rho}$ and $D,\alpha>0$ such that
	\begin{align}
		& E(w) =E_*-D, & & 0<D <\delta,\label{eq:gluing-energy}\\
		& P(w) =-\alpha, & & 0<\alpha \leq6D,\label{eq:gluing-P}
	\end{align}
	and
	\begin{equation}\label{eq:gluing-Q}
		\mathcal Q(w)\leq\sum_{j=1}^J\mathcal Q(U^j)+4L_0(b)\leq q<0.
	\end{equation}
\end{lemma}

\begin{proof}
First, we construct a pair $T_L$ with component masses $b_i$.
Fix a nonnegative $\eta\in C_c^\infty(\R)$ with
$\|\eta\|_{L^2(\R)}=1$. For $L>0$, define
$T_L=(T_{1,L},T_{2,L})$ by
\begin{equation}\label{eq:long-wave-tail}
    T_{i,L}(y,z):=\sqrt{b_i}\,\phi_0(y)L^{-\frac{1}{2}}\eta(\frac{z}{L}).
\end{equation}
The normalization of $\eta$ and $\phi_0$ gives
$\|T_{i,L}\|_2^2=b_i$. Set
\[
A_b:=(b_1+b_2)\|\eta'\|_{L^2(\R)}^2,
\qquad
c_b:=\frac{\|\eta\|_{L^4(\R)}^4}{2\pi}
\left(\frac{\mu_1}{4}b_1^2+\frac{\mu_2}{4}b_2^2+\beta b_1b_2\right).
\]
Using $\|\phi_0\|_{L^4(\R^2)}^4=\frac{1}{2\pi}$ and
$\|\nabla_y\phi_0\|_{L^2(\R^2)}^2=\||y|\phi_0\|_{L^2(\R^2)}^2=1$, we obtain
\[
\begin{aligned}
    & E(T_L)=L_0(b)-\frac{c_b}{L}+\frac{A_b}{2L^2},\\
    & P(T_L)=-\frac{3c_b}{L}+\frac{A_b}{L^2},\\
    & \mathcal Q(T_L)=4L_0(b)-\frac{9c_b}{L}+\frac{2A_b}{L^2}.
\end{aligned}
\]
If $b\neq0$, then $c_b>0$. Define
$D_{\rm tail}:=L_0(b)-E(T_L)$ and $\alpha_{\rm tail}:=-P(T_L)$.
By choosing $L>\frac{A_b}{c_b}$ sufficiently large, we obtain
\begin{equation}\label{eq:tail-budget}
    0<D_{\rm tail}<\frac{\delta}{2},\qquad
    0<\alpha_{\rm tail}\leq6D_{\rm tail},\qquad
    \mathcal Q(T_L)\leq4L_0(b).
\end{equation}
When $b=0$, set $T_L=0$ and
$D_{\rm tail}=\alpha_{\rm tail}=0$.

Next, we combine $T_L$ with the profiles $U^1,\ldots,U^J$.
Choose $\delta_1,\ldots,\delta_J>0$ with
$\sum_{j=1}^J\delta_j<\delta-D_{\rm tail}$ and set $w^0=T_L$.
For $j=1,\ldots,J$, define
\[
w_i^j:=\bigl((w_i^{j-1})^2+(\tau_{s_j}U_i^j)^2\bigr)^{\frac{1}{2}},
\qquad i=1,2,
\]
where $s_j\in\R$ will be chosen below.
By \cref{lem:two-block-quadrature-separation}, the quantities
\[
\begin{aligned}
    & \kappa_j:=K(w^{j-1})+K(U^j)-K(w^j),\\
    & \gamma_j:=B(w^j)-B(w^{j-1})-B(U^j)
\end{aligned}
\]
are nonnegative. Set $d_j:=\frac{\kappa_j}{2}+\gamma_j$.
At each stage, $w^{j-1}$ is fixed. By the limits in
\eqref{eq:two-block-separation-limit}, we can choose $s_j$ so that
$d_j<\delta_j$. When $b=0$, we take $s_1=0$, for which
$w^1=U^1$ and $d_1=0$.

Finally, let $w=w^J$. By the mass and potential identities in
\cref{lem:two-block-quadrature-separation}, we obtain
\[
\|w_i\|_2^2=b_i+\sum_{j=1}^J\|U_i^j\|_2^2=\eps\rho_i,
\qquad
W(w)=W(T_L)+\sum_{j=1}^J W(U^j).
\]
In particular, $w\in\cS_{\eps\rho}$.
To obtain the required estimates, define
\begin{align*}
    & \Delta_K:=\sum_{j=1}^J\kappa_j, && \Delta_B:=\sum_{j=1}^J\gamma_j,\\
    & D_{\rm int}:=\frac{1}{2}\Delta_K+\Delta_B,
    && \alpha_{\rm int}:=\Delta_K+3\Delta_B.
\end{align*}
These definitions give $0\leq\alpha_{\rm int}\leq3D_{\rm int}$ and
$D_{\rm int}=\sum_{j=1}^Jd_j<\delta-D_{\rm tail}$.
Since $P(U^j)=0$, summing the defining identities for
$\kappa_j$ and $\gamma_j$ gives
\[
\begin{aligned}
    & E(w)=E_*-D_{\rm tail}-D_{\rm int},\\
    & P(w)=-\alpha_{\rm tail}-\alpha_{\rm int},\\
    & \mathcal Q(w)=\sum_{j=1}^J\mathcal Q(U^j)+\mathcal Q(T_L)
      -2\Delta_K-9\Delta_B.
\end{aligned}
\]
Set $D=D_{\rm tail}+D_{\rm int}$ and
$\alpha=\alpha_{\rm tail}+\alpha_{\rm int}$.
It remains to verify that $D>0$ and $\alpha>0$.
If $b\neq0$, then $D_{\rm tail},\alpha_{\rm tail}>0$ by
\eqref{eq:tail-budget}. If $b=0$, the assumption gives $J\geq2$.
In this case, choose $i,k\in\{1,2\}$ such that $U_i^1\not\equiv0$ and
$U_k^2\not\equiv0$. By \cref{lem:profile-equations-pohozaev}, both
components are strictly positive. Therefore,
\[
\int_{\R^3}(U_i^1)^2(\tau_{s_2}U_k^2)^2\dd x>0
\qquad(s_2\in\R).
\]
The coefficient of this integral in $\gamma_2$ is $\frac{\mu_i}{2}$ if
$i=k$ and $\beta$ if $i\neq k$. Thus $\gamma_2>0$, which gives
$D>0$ and $\alpha>0$ also in this case.
The preceding bounds now yield
\[
0<D<\delta,\qquad
0<\alpha\leq6D_{\rm tail}+3D_{\rm int}\leq6D.
\]
Moreover, \eqref{eq:tail-budget} and \cref{lem:residual-ledger} give
\[
\mathcal Q(w)\leq\sum_{j=1}^J\mathcal Q(U^j)+4L_0(b)\leq q<0.
\]
Together with the identities for $E(w)$ and $P(w)$, these estimates
prove \eqref{eq:gluing-energy}--\eqref{eq:gluing-Q}.
\end{proof}

The pair $w$ constructed in \cref{lem:finite-quadrature-gluing}
satisfies $P(w)=-\alpha<0$. To obtain an element of
$\mathcal P_{\eps\rho}^-$, we use the following lemma, which also
bounds the energy increase by $C\alpha^2$.

\begin{lemma}
	\label{lem:negative-sheet-projection}
	Fix $\eps,\delta_0>0$, and let $\mathcal W\subset\cS_{\eps\rho}$ be
	bounded in $\cH^2$. There exist constants $\alpha_*,C_{\mathrm{pr}}>0$,
	depending only on $\mathcal W$ and $\delta_0$, such that for every
	$w\in\mathcal W$ with
	\[
	P(w)=-\alpha,\qquad 0<\alpha<\alpha_*,\qquad
	\mathcal Q(w)\leq-\delta_0,
	\]
	there exists
	$t_\alpha\in[1-C_{\mathrm{pr}}\alpha,1)$ such that
	\begin{equation}\label{eq:projected-negative-sheet}
		P(t_\alpha\star w)=0,
		\qquad
		\mathcal Q(t_\alpha\star w)\leq-\frac{\delta_0}{2},
	\end{equation}
	and
	\begin{equation}\label{eq:projection-energy-cost}
		0\leq E(t_\alpha\star w)-E(w)
		\leq C_{\mathrm{pr}}\alpha^2.
	\end{equation}
\end{lemma}

\begin{proof}
Let $p(t)=P(t\star w)$ for $t>0$. Then
\[
p(1)=-\alpha,\qquad p'(t)=\frac{\mathcal Q(t\star w)}{t}.
\]
Since $\mathcal W$ is bounded in $\cH^2$, the quantities $K(w)$,
$W(w)$, and $B(w)$ are uniformly bounded for $w\in\mathcal W$.
By the scaling formula for $\mathcal Q$ and the mean value theorem,
we obtain $C_{\mathcal W}>0$ such that
\[
\bigl|\mathcal Q(t\star w)-\mathcal Q(w)\bigr|
\leq C_{\mathcal W}|t-1|,
\qquad w\in\mathcal W,\quad t\in[\frac{1}{2},1].
\]
Choose $r_0>0$ depending only on $\mathcal W$ and $\delta_0$ such that
\[
r_0<\min\left\{\frac{1}{2},\frac{\delta_0}{2C_{\mathcal W}}\right\}.
\]
Then, for every $w\in\mathcal W$ under the hypotheses of the lemma,
we obtain
\[
\mathcal Q(t\star w)
\leq\mathcal Q(w)+C_{\mathcal W}|t-1|
\leq-\frac{\delta_0}{2}
\qquad(t\in[1-r_0,1]).
\]
Hence
\[
p'(t)\leq-\frac{\delta_0}{2}
\qquad(t\in[1-r_0,1]).
\]

Next, choose $0<\alpha_*\leq\frac{r_0\delta_0}{2}$ so that
\[
1-\frac{2\alpha}{\delta_0}\in[1-r_0,1)
\qquad\text{whenever }0<\alpha<\alpha_*.
\]
Consequently, the derivative estimate above gives
\[
p\left(1-\frac{2\alpha}{\delta_0}\right)
=p(1)-\int_{1-\frac{2\alpha}{\delta_0}}^1p'(t)\dd t
\geq-\alpha+\alpha=0.
\]
Together with $p(1)=-\alpha<0$, this inequality allows us to apply
the intermediate value theorem and obtain
\[
t_\alpha\in[1-\frac{2\alpha}{\delta_0},1),
\qquad p(t_\alpha)=0.
\]
By the bound for $\mathcal Q$ on $[1-r_0,1]$, we also obtain
$\mathcal Q(t_\alpha\star w)\leq-\frac{\delta_0}{2}$.

Finally, we estimate the energy increase caused by this scaling.
By the monotonicity of $p$, we have
\[
-\alpha\leq p(t)\leq0\qquad(t\in[t_\alpha,1]).
\]
Using $t_\alpha\geq1-r_0>\frac{1}{2}$ and the bound for
$1-t_\alpha$, we obtain
\[
0\leq E(t_\alpha\star w)-E(w)
=-\int_{t_\alpha}^1\frac{p(t)}{t}\dd t
\leq2\alpha(1-t_\alpha)\leq\frac{4}{\delta_0}\alpha^2.
\]
Set $C_{\mathrm{pr}}=\frac{4}{\delta_0}$.
Then $C_{\mathrm{pr}}\alpha_*\leq2r_0<1$ and
$1-C_{\mathrm{pr}}\alpha\leq1-\frac{2\alpha}{\delta_0}$.
Thus $t_\alpha\in[1-C_{\mathrm{pr}}\alpha,1)\subset(0,\infty)$,
which completes the proof of \eqref{eq:projected-negative-sheet}
and \eqref{eq:projection-energy-cost}.
\end{proof}

Now we apply \cref{lem:finite-quadrature-gluing,lem:negative-sheet-projection}
to prove that $d_\eps$ is attained.

\begin{theorem}
	\label{thm:high-level-attainment}
	Let $\mu_1,\mu_2,\beta>0$ and $\rho\in(0,\infty)^2$.
	For every sufficiently small $\eps>0$, the level $d_\eps$
	is attained by a positive pair $u_\eps^-\in\mathcal P_{\eps\rho}^-$.
	Moreover, every minimizer of $d_\eps$ on $\mathcal P_{\eps\rho}^-$
	is a constrained critical point of $E$ on $\cS_{\eps\rho}$.
\end{theorem}

\begin{proof}
	Let $\{u_n\}$ be the sequence in \cref{prop:high-PS}, and let
	$U^j,b,\ell,q$ be as in
	\cref{lem:axial-profile-decomposition,lem:residual-ledger}.

	\smallskip
	\noindent\textbf{Step 1.}
	First, we prove that $J=1$ and $b=0$.
	Suppose, by contradiction, that $b\neq0$ or $J\geq2$.
	For $0<\delta<1$, let $w$ be the pair constructed in
	\cref{lem:finite-quadrature-gluing}.
	By \eqref{eq:gluing-energy} and \eqref{eq:gluing-P}, we obtain
	\[
	K(w)+5W(w)=6E(w)-2P(w)=6E_*-6D+2\alpha
	\leq6E_*+6,
	\]
	where $0<D<\delta<1$ and $0<\alpha\leq6D$.
	Since $W(w)\geq0$ and $\|w_i\|_2^2=\eps\rho_i$ for $i\in\{1,2\}$,
	we obtain
	\[
	\|w\|_{\cH^2}^2
	=K(w)+W(w)+\eps(\rho_1+\rho_2)
	\leq6E_*+6+\eps(\rho_1+\rho_2).
	\]
	This bound is independent of $\delta$.
	Therefore, we can choose the constants $\alpha_*,C_{\mathrm{pr}}$ in
	\cref{lem:negative-sheet-projection} independently of $\delta$,
	with $\delta_0=-q>0$.
	Now choose
	\[
	0<\delta<\min\left\{1,\frac{\alpha_*}{6},
	\frac{1}{36C_{\mathrm{pr}}}\right\}
	\]
	and apply \cref{lem:finite-quadrature-gluing} to obtain
	$0<D<\delta$ and $0<\alpha\leq6D<\alpha_*$.
	By \eqref{eq:gluing-Q}, we also have
	$\mathcal Q(w)\leq q=-\delta_0$.
	Hence \cref{lem:negative-sheet-projection} gives
	$t_\alpha\star w\in\mathcal P_{\eps\rho}^-$ such that
	\[
	\begin{aligned}
	&E(t_\alpha\star w)\leq E_*-D+C_{\mathrm{pr}}\alpha^2
	\leq E_*-D+36C_{\mathrm{pr}}D^2<E_*,\\
	&E_*\leq\sum_{j=1}^J E(U^j)+\ell=d_\eps.
	\end{aligned}
	\]
	By the definition of $d_\eps$, this would imply
	$d_\eps\leq E(t_\alpha\star w)<d_\eps$.
	Since $J\geq1$, this contradiction proves that $J=1$ and $b=0$.

	\smallskip
	\noindent\textbf{Step 2.}
	Next, we prove that $\ell=0$ and establish strong convergence to a
	positive minimizer of $d_\eps$.
	Set $U=U^1$. Since $b=0$, \eqref{eq:missing-profile-mass} gives
	$\|U_i\|_2^2=\eps\rho_i$ for $i\in\{1,2\}$.
	Suppose, by contradiction, that $\ell>0$.
	By \eqref{eq:Q-ledger} and \eqref{eq:ell-profile}, we obtain
	\[
	\mathcal Q(U)=q-4\ell<0,
	\qquad E(U)=d_\eps-\ell<d_\eps.
	\]
	Together with $P(U)=0$ and the mass constraints, these inequalities
	give $U\in\mathcal P_{\eps\rho}^-$ and
	$d_\eps\leq E(U)<d_\eps$, a contradiction.
	As $\ell\geq L_0(0)=0$ by \cref{lem:residual-ledger}, we conclude that
	$\ell=0$. Consequently, $\mathcal Q(U)=q<0$ and $E(U)=d_\eps$,
	which proves that $U$ is a minimizer of $d_\eps$ on
	$\mathcal P_{\eps\rho}^-$.

	Now the decomposition reads $u_n=\tau_{z_n^1}U+r_n$.
	From \cref{lem:residual-ledger} with $b=0$ and $\ell=0$, we obtain
	\[
	\|r_n\|_{\cH^2}^2
	=K(r_n)+W(r_n)+\sum_i\|r_{i,n}\|_2^2\longrightarrow0.
	\]
	Since axial translations preserve the $\cH^2$ norm, it follows that
	\[
	\tau_{-z_n^1}u_n\longrightarrow U\quad\text{strongly in }\cH^2.
	\]
	Moreover, \eqref{eq:profile-limit-system} shows that $U$ is a
	constrained critical point of $E$ on $\cS_{\eps\rho}$.
	Since both masses are positive, \cref{lem:profile-equations-pohozaev}
	implies that both components of $U$ are strictly positive.
	Thus $u_\eps^-:=U$ is the required positive normalized solution.

	\smallskip
	\noindent\textbf{Step 3.}
	Finally, we show that every minimizer of $d_\eps$ on
	$\mathcal P_{\eps\rho}^-$ is a constrained critical point of $E$
	on $\cS_{\eps\rho}$.
	Let $\overline U\in\mathcal P_{\eps\rho}^-$ attain $d_\eps$.
	Recall the augmented functionals in
	\eqref{eq:augmented-energy}--\eqref{eq:augmented-PQ}.
	At $(1,\overline U)$, we have
	\[
	\partial_t\widetilde E(1,\overline U)=P(\overline U)=0,\qquad
	\partial_t\widetilde P(1,\overline U)=\mathcal Q(\overline U)<0.
	\]
	The mass constraints do not depend on $t$, and their derivatives
	are linearly independent. Together with
	$\partial_t\widetilde P(1,\overline U)\neq0$, this shows that the
	derivatives of all three constraints are linearly independent.
	By continuity, $\widetilde{\mathcal Q}(t,w)<0$ near $(1,\overline U)$.
	Thus every nearby $(t,w)$ satisfying the mass constraints and
	$\widetilde P(t,w)=0$ satisfies $t\star w\in\mathcal P_{\eps\rho}^-$.
	For these points, the definition of $d_\eps$ gives
	\[
	\widetilde E(t,w)=E(t\star w)
	\geq d_\eps=\widetilde E(1,\overline U).
	\]
	Consequently, $(1,\overline U)$ is a local minimizer subject to
	the three constraints.
	By the Lagrange multiplier rule, there exist
	$\eta,\lambda_1,\lambda_2\in\R$, where $\eta$ corresponds to
	$\widetilde P=0$.
	The derivative with respect to $t$ gives
	$0=-\eta\mathcal Q(\overline U)$.
	Since $\mathcal Q(\overline U)<0$, we obtain $\eta=0$ and hence
	\[
	E'(\overline U)=\sum_{i=1}^2\lambda_i\overline U_i e_i
	\quad\text{in }(\cH^2)'.
	\]
	This proves that $\overline U$ is a constrained critical point of
	$E$ on $\cS_{\eps\rho}$.
	If $\lambda_i$ and $\widehat\lambda_i$ are two choices of multipliers,
	subtracting their equations and testing with $\overline U_i e_i$ gives
	$(\lambda_i-\widehat\lambda_i)\eps\rho_i=0$.
	Since $\eps\rho_i>0$, each multiplier is unique.
\end{proof}

\begin{proposition}\label{prop:high-asymptotics}
	Every family of positive minimizers $u_\eps^-$ of $d_\eps$ on
	$\mathcal P_{\eps\rho}^-$ satisfies the limits in
	\eqref{eq:main-high-limits} as $\eps\to0$.
\end{proposition}

\begin{proof}
	First, we prove that $\eps d_\eps\to m_0(\rho)$ as $\eps\to0$.
	Fix $\delta>0$. As in the proof of \cref{lem:fold-gap}, choose a smooth
	compactly supported $v\in\cS_\rho^0$ with $P_0(v)=0$ and
	$E_0(v)<m_0(\rho)+\delta$.
	For this fixed $v$ and all sufficiently small $\eps>0$,
	\eqref{eq:trial-energy-bound} and the lower bound
	in \eqref{eq:d-upper-rough} give
	\[
	\begin{aligned}
	&\eps d_\eps\leq E_0(v)+\frac{50}{81}\eps^4W(v),\\
	&m_0(\rho)\leq\liminf_{\eps\to0}\eps d_\eps
	\leq\limsup_{\eps\to0}\eps d_\eps\leq m_0(\rho)+\delta.
	\end{aligned}
	\]
	Letting $\delta\downarrow0$ yields $\eps d_\eps\to m_0(\rho)$.

	Next, we show that
	\[
	r_\eps:=\frac{W(u_\eps^-)}{K(u_\eps^-)}\longrightarrow0
	\qquad\text{as }\eps\to0.
	\]
	Since $u_\eps^-\in\mathcal P_{\eps\rho}^-$, we have
	$0\leq r_\eps<\frac{1}{5}$. By \eqref{eq:fiber-lower-f}, we obtain
	\[
	1\leq f(r_\eps)\leq\frac{\eps d_\eps}{m_0(\rho)}\to1,
	\qquad f(r)=(1-r)^2(1+5r).
	\]
	Since $f$ is strictly increasing on $[0,\frac{1}{5}]$ and $f(0)=1$,
	we conclude that $r_\eps\to0$.

	Finally, the Pohozaev identity
	$K(u_\eps^-)+5W(u_\eps^-)=6d_\eps$ in \eqref{eq:sheet-identities}
	and the preceding limits give
	\[
	\eps K(u_\eps^-)=\frac{6\eps d_\eps}{1+5r_\eps}\to6m_0(\rho),
	\qquad
	\eps W(u_\eps^-)=r_\eps\eps K(u_\eps^-)\to0.\qedhere
	\]
\end{proof}

Now we combine the mountain pass characterization with the preceding
results to prove \cref{thm:main}.

\begin{proof}[Proof of \cref{thm:main}]
	Take $\eps_0>0$ below the thresholds in
	\cref{prop:mp-level-identification,thm:high-level-attainment}.
	Combining these results, we obtain a positive constrained critical
	point $u_\eps^-\in\mathcal P_{\eps\rho}^-$ satisfying
	\[
	E(u_\eps^-)=m_\eps=d_\eps>
	\max\{E(p_\eps),E(e_\eps)\}
	\qquad \eps\in(0,\eps_0).
	\]
	Finally, \cref{prop:high-asymptotics} gives the limits in
	\eqref{eq:main-high-limits}.
	The choice of $\eps_0$ depends only on the parameters in the statement.
\end{proof}

Finally, we use the compactness argument above to prove
\cref{cor:negative-sheet-compactness} for arbitrary minimizing sequences.

\begin{proof}[Proof of \cref{cor:negative-sheet-compactness}]
	First, let $\{v_n\}\subset\mathcal P_{\eps\rho}^-$ be a minimizing sequence
	for $d_\eps$. Applying the construction of \cref{prop:high-PS}
	to $\{|v_n|\}$ gives $u_n=t_n\star w_n$, with $t_n\to1$ and
	$\|w_n-|v_n|\|_{\cH^2}\to0$.
	By the compactness argument in the proof of
	\cref{thm:high-level-attainment}, up to a subsequence there exist
	$\{a_n\}\subset\R$ and a positive minimizer
	$U\in\mathcal P_{\eps\rho}^-$ such that
	$\tau_{-a_n}u_n\to U$ in $\cH^2$.

	Next, we obtain convergence of $\tau_{-b_n}|v_n|$, where $b_n=t_na_n$.
	The identity
	$\tau_{-b_n}w_n=t_n^{-1}\star(\tau_{-a_n}u_n)$ and the joint
	continuity of scaling imply $\tau_{-b_n}w_n\to U$ in $\cH^2$.
	Since translations preserve the norm and $\|w_n-|v_n|\|_{\cH^2}\to0$,
	we obtain $\tau_{-b_n}|v_n|\to U$ in $\cH^2$.

	Finally, we prove strong convergence of $\tau_{-b_n}v_n$.
	By boundedness and the Rellich-Kondrachov theorem, up to a subsequence
	we have $\tau_{-b_n}v_n\rightharpoonup V$ in $\cH^2$ and
	$\tau_{-b_n}v_n\to V$ in $L^2_{\mathrm{loc}}$.
	Comparing these limits with $\tau_{-b_n}|v_n|\to U$ gives $|V|=U$.
	Since $U_i\in C^2_{\mathrm{loc}}(\R^3)$ is strictly positive,
	$1/U_i\in W_{\mathrm{loc}}^{1,\infty}(\R^3)$. Hence $V_i/U_i\in H^1_{\mathrm{loc}}(\R^3)$.
	Applying the chain rule for weak derivatives to $(V_i/U_i)^2=1$ gives
	$2(V_i/U_i)\nabla(V_i/U_i)=0$ almost everywhere.
	Therefore, $\nabla(V_i/U_i)=0$, and connectedness of $\R^3$
	gives $V_i=\pm U_i$.
	Moreover,
	\[
	\|\tau_{-b_n}v_n\|_{\cH^2}
	=\|\tau_{-b_n}|v_n|\|_{\cH^2}
	\longrightarrow\|U\|_{\cH^2}=\|V\|_{\cH^2}.
	\]
	Together with weak convergence, this proves $\tau_{-b_n}v_n\to V$
	in $\cH^2$. Since $|V|=U$, we have
	$V\in\mathcal P_{\eps\rho}^-$ and $E(V)=d_\eps$.
	If $v_n\geq0$ for every $n$, then $V=U>0$.
\end{proof}

\section{Conclusion}

In this paper, we studied positive normalized solutions of a cooperative cubic Schr\"dinger system with partial confinement and two prescribed component masses. For fixed positive interaction coefficients and a fixed positive mass vector, we proved that the system admits a positive normalized solution at a mountain pass level when the masses are sufficiently small. We also identified this level with the infimum of the energy on the negative part of the Pohozaev set. In addition to the normalized ground state obtained previously by Jia et al. \cite{JiaLiLuo2020}, we established the existence of a second positive normalized solution at the mountain pass level.

A key part of the proof is the compactness analysis for minimizing sequences under the two mass constraints. The axial translation invariance allows sequences to split into several profiles, while mass and quadratic energy may remain in the remainder. We address these possibilities by combining an axial profile decomposition with a construction that combines profiles and restores any missing mass. This yields strong compactness of every minimizing sequence for the negative Pohozaev level, up to a subsequence and common axial translations. We also obtain the small-mass limits
\[
\varepsilon E(u_{\varepsilon}^-)\to m_0(\rho),~
\varepsilon K(u_{\varepsilon}^-)\to 6m_0(\rho),~
\varepsilon W(u_{\varepsilon}^-)\to 0
\quad\text{as }\varepsilon\to0.
\]

An interesting direction for future work is to extend these results to coupled Schr\"odinger systems with general mass supercritical and Sobolev subcritical nonlinearities. When the nonlinear terms have different scaling exponents, the mountain pass geometry, compactness estimates, and the small mass energy scale require further analysis.\\

\noindent\textbf{Author Contributions}
Wan-Yue Tang: Writing---original draft, Validation, Formal analysis, Investigation, Conceptualization.
Hong-Rui Sun: Writing---review and editing, Validation, Supervision, Methodology, Funding acquisition, Conceptualization.

\medskip
\noindent\textbf{Funding}
H.~R.~Sun was partially supported by the National Natural Science Foundation of China (Grant No. 11671181) and the Natural Science Foundation of Gansu Province of China (Grant Nos. 21JR7RA535 and 24JRRA414).

\medskip
\noindent\textbf{Data Availability}
No datasets were generated or analysed during the current study.

\medskip
\noindent\textbf{Materials Availability}
Not applicable.

\section*{Declarations}

\noindent\textbf{Conflict of interest}
The authors declare no conflict of interest.

\medskip
\noindent\textbf{Ethics approval}
Not applicable.

\medskip
\noindent\textbf{Consent for publication}
All authors have reviewed and approved the final version of the manuscript for submission and publication in the journal.

\end{document}